\documentclass[12pt]{article}
\usepackage{amssymb}
\usepackage{amsmath}
\usepackage{latexsym}

\numberwithin{equation}{section}
\newtheorem{thm}[equation]{Theorem}

\newtheorem{defn}[equation]{Definition}
\newtheorem{rem}[equation]{Remark}
\newtheorem{lem}[equation]{Lemma}
\newtheorem{corol}[equation]{Corollary}
\newtheorem{example}[equation]{Example}
\title{Spectral stability of higher order uniformly elliptic operators\footnote{This paper has been published in: {\it Sobolev spaces in mathematics}. II,  Int. Math. Ser. (N. Y.), {\bf 9}, 69-102, Springer, New York, 2009 }}
\author{Victor I. Burenkov and Pier Domenico Lamberti}

\date{\ }

\begin{document}

\newcommand{\rea}{\mathbb{R}}

\maketitle

%
%
%

\noindent
{\bf Abstract:}
We prove estimates for the variation of the eigenvalues of uniformly elliptic operators with homogeneous Dirichlet or Neumann boundary conditions
 upon variation of the open set on which an operator is defined. We consider operators of arbitrary even order and open sets
 admitting arbitrary strong degeneration.  The main estimate
 is expressed via a natural and easily computable distance between open sets with continuous boundaries.
 Another estimate
 is obtained via the lower Hausdorff-Pompeiu deviation of the boundaries, which in general may be much smaller than the usual Hausdorff-Pompeiu distance.
Finally, in the case of diffeomorphic open sets we obtain an estimate even without the assumption of continuity of the boundaries.
\\

\vspace{11pt}


\vspace{6pt}

\section{Introduction}

In this paper we consider the eigenvalue problem for the operator

\begin{equation}
\label{classic}
Hu=
 (-1)^m\sum_{|\alpha |=|\beta |=m} D^{\alpha }\left(A_{\alpha \beta}(x)D^{\beta }u  \right),\
\ \ x\in\Omega ,
\end{equation}
subject to homogeneous Dirichlet or Neumann  boundary conditions, where $m\in {\mathbb{N}}$, $\Omega $ is a bounded  open set in ${\mathbb{R}}^N$ and the coefficients $A_{\alpha \beta }$ are Lipschitz continuous  functions satisfying the uniform ellipticity condition (\ref{elp}) on $\Omega$. For a precise statement of the eigenvalue problem, see  Definition \ref{defset} and Theorem \ref{setting}.

We consider open sets $\Omega$ for which the spectrum is discrete and can be represented by means of a non-decreasing sequence
of non-negative eigenvalues
$$
\lambda_1[\Omega ]\le \lambda_2[\Omega ]\le \dots \le \lambda_n[\Omega]\le \dots
$$
where each eigenvalue is repeated as many times as its multiplicity.

In this paper we prove estimates for the variation
$$
|\lambda_n[\Omega_1]-\lambda_n[\Omega_2]|
$$
of the eigenvalues corresponding to two open sets $\Omega_1$, $\Omega_2$.

There is  vast literature on spectral stability problems for elliptic operators (see {\it e.g.,} Hale~\cite{hale}, Henry~\cite{he} for references). However, little attention has
been devoted to the problem of spectral stability for higher order operators and  in particular  to the problem of finding  explicit
qualified estimates for the variation of the eigenvalues.  Moreover,
most of the existing qualified estimates for second order operators were obtained under certain  regularity assumptions on the boundaries.

Our analysis comprehends operators of arbitrary even order, with  homogeneous Dirichlet or  Neumann  boundary conditions, and open sets
admitting arbitrarily strong degeneration. In fact, we consider bounded open sets  whose boundaries are just locally the subgraphs of  continuous functions. We only require that the `atlas'  ${\mathcal{A}}$, with the help of which such open sets are described, is fixed:  we denote by $C({\mathcal{A}})$ the family of all such  open sets (see Definition \ref{class}). In $C({\mathcal{A}})$ we introduce a
natural metric $d_{{\mathcal{A}}}$ (the `atlas distance') which can be easily computed. Given two open sets $\Omega_1, \Omega_2\in C({\mathcal{A}})$, the distance
  $d_{{\mathcal{A}}}(\Omega_1,\Omega_2)$ is just the maximum of the sup-norms of  the differences of the functions describing locally the boundaries of
  $\Omega_1$ and $\Omega_2$ (see Definition \ref{dev}).

  The first main result of the paper is that for both  Dirichlet and Neumann boundary conditions the eigenvalues of (\ref{classic}) are locally Lipschitz continuous
  functions of the open set $\Omega \in C({\mathcal{A}})$ with respect to the atlas distance  $d_{{\mathcal{A}}}$. Namely, in Theorems \ref{dirthm} and \ref{neuthm} we prove that
  for each $n\in {\mathbb{N}}$ there exist $c_n, \epsilon_n>0$ such that for both Dirichlet and Neumann boundary conditions the estimate
\begin{equation}
\label{mainintr}
  |\lambda_n[\Omega_1]-\lambda_n[\Omega_2]|\le c_n d_{{\mathcal{A}}}(\Omega_1,\Omega_2)
\end{equation}
holds for all open sets $\Omega_1, \Omega_2\in C({\mathcal{A}})$ satisfying $d_{{\mathcal{A}}}(\Omega_1,\Omega_2)<\epsilon_n$.

By estimate (\ref{mainintr}) we deduce an estimate expressed in terms of the lower Hausdorff-Pompeiu deviation of the boundaries
$$
d_{{\mathcal{H}}{\mathcal{P}}}(\partial \Omega_1 , \partial \Omega_2 )=\min \left\{\sup_{x\in \partial \Omega_1}d(x, \partial \Omega_2), \sup_{x\in \partial \Omega_2}d(x, \partial \Omega_1)  \right\}.
$$

To do so, we restrict our attention to smaller families of open sets in  $C({\mathcal{A}})$. Namely, for a fixed $M>0$ and $\omega :[0,\infty [\to [0,\infty [$ satisfying very weak natural conditions, we consider those open sets $\Omega $ in $C({\mathcal{A}})$ for which   any of the functions $\bar x\mapsto g(\bar x)$, describing locally the boundary
of $\Omega$, satisfies the condition
$$|g(\bar x )-g(\bar y )|\le M \omega (|\bar x -\bar y  |),$$
for all appropriate $\bar x, \bar y$:
we denote by $C_M^{\omega (\cdot )}({\mathcal{A}})$ the family of all such open sets (see Definition \ref{omegafam}). For instance, if $0<\alpha \le 1$ and $\omega (t)=t^{\alpha }$ for all $t\geq 0$
 then we obtain open sets with H\"{o}lder continuous boundaries of exponent $\alpha$: this class is denoted below by $C_M^{0, \alpha }({\mathcal{A}})$.   It is possible to  choose  a function $\omega $ going to zero arbitrarily slowly which allows dealing with open sets with  arbitrarily sharp cusps.

The second main result of the paper is for open sets $\Omega_1,\Omega_2\in C_M^{\omega (\cdot )}({\mathcal{A}})$. Namely, in Theorem \ref{maincor} we prove that for each $n\in {\mathbb{N}}$ there exist $c_n, \epsilon_n>0$ such that
\begin{equation}
\label{mainintr1}
  |\lambda_n[\Omega_1]-\lambda_n[\Omega_2]|\le c_n \omega (d_{{\mathcal{H}} {\mathcal{P}}}(\partial \Omega_1,\partial \Omega_2)),
\end{equation}
for all open sets $\Omega_1, \Omega_2\in  C_M^{\omega (\cdot )}({\mathcal{A}})$ satisfying $d_{{\mathcal{H}} {\mathcal{P}}}(\partial \Omega_1,\partial \Omega_2)<\epsilon_n$.

In particular in Corollary \ref{lastcor} we deduce  that if $\Omega _1, \Omega_2\in C_M^{\omega (\cdot )}({\mathcal{A}})$ satisfy
$$
(\Omega_1)_{\epsilon }\subset \Omega_2 \subset (\Omega_1)^{\epsilon }\ \ \ \ {\rm or }\ \ \ \ (\Omega_2)_{\epsilon }\subset \Omega_1 \subset (\Omega_2)^{\epsilon },
$$
where
 $\epsilon >0 $ is sufficiently small
then
\begin{equation}
\label{mainintr2}
  |\lambda_n[\Omega_1]-\lambda_n[\Omega_2]|\le c_n \omega (\epsilon ).
\end{equation}
Here  $\Omega _{\epsilon }=\left\{x\in \Omega :\ d(x,\partial \Omega )>\epsilon   \right\}$,
$\Omega ^{\epsilon }=\left\{x\in {\mathbb{R}}^N :\ d(x, \Omega )<\epsilon   \right\}$, for any set $\Omega$ in ${\mathbb{R}}^N$.

In the case $\Omega_1, \Omega_2\in C^{0,\alpha }_M({\mathcal{A}})$ estimate (\ref{mainintr2}) takes the form
\begin{equation}
\label{mainintr3}
  |\lambda_n[\Omega_1]-\lambda_n[\Omega_2]|\le c_n \epsilon^{\alpha } .
\end{equation}

In the case of Dirichlet boundary conditions and $m=1$ some estimates of the form (\ref{mainintr3}) were obtained in Davies~\cite{da2000}
under the assumption that a certain Hardy-type inequality is satisfied on $\Omega_1$ (see also Pang~\cite{pang1}).
In the case of Dirichlet  boundary conditions and $m=1$,
estimate (\ref{mainintr3}) was  proved in \cite{buladir}.
In the case of Dirichlet boundary conditions, $m=2$ and  open sets with sufficiently smooth boundaries an estimate of the form   (\ref{mainintr3}) was obtained
in Barbatis~\cite{bar}.

In the case of Neumann boundary conditions and $m=1$, estimate (\ref{mainintr3}) was proved in  Burenkov and Davies~\cite{buda} for open sets $\Omega _1, \Omega _2\in C^{0,\alpha }_M({\mathcal{A}})$ satisfying $(\Omega_1)_{\epsilon }\subset \Omega_2 \subset \Omega_1 $.
We remark that the result in \cite{buda} concerns only inner deformations of an open set and second order elliptic operators. Moreover, the proof
in \cite{buda} is based on the ultracontractivity which  holds for second order elliptic operators in open sets with H\"{o}lder continuous boundaries. Since  ultracontractivity
is not guaranteed for more general  open sets, we had to develop a different method.

The third main result of the paper
concerns  the case  $\Omega_1=\Omega $ and $\Omega_2=\phi (\Omega )$,  where $\phi$ is  a suitable diffeomorphism
of class $C^m$. In this case we make  very weak  assumptions on $\Omega$: if $m=1$ it is just the requirement that $H$ has discrete spectrum.
Under such general assumptions we prove that for both Dirichlet and Neumann boundary conditions there exists a constant
$c>0$ independent of $n$ such that
$$
|\lambda_{n}[\Omega ]-\lambda_{n}[\phi (\Omega )]|\le c(1+\lambda_{n}[\Omega ]) \max_{0\le |\alpha |\le m}   \| D^{\alpha }(\phi -{\rm Id} )\|_{L^{\infty }(\Omega )},
$$
if $\max_{0\le |\alpha |\le m}   \| D^{\alpha }(\phi -{\rm Id} )\|_{L^{\infty }(\Omega )}<c^{-1}$ (see Theorem~\ref{diffeothm} and Corollary~\ref{diffeocorol}).

The paper is organized as follows: in Section~\ref{section2} we introduce some notation and we formulate the eigenvalue problem for operator (\ref{classic}); in Section~\ref{section3} we define the class of open sets under consideration;  in Section~\ref{section4} we consider the case of diffeormorphic open sets; in Section~\ref{section5} we prove estimate (\ref{mainintr}) for  Dirichlet boundary conditions; in Section~\ref{section6} we prove estimate (\ref{mainintr}) for  Neumann boundary conditions;
in Section~\ref{section7} we prove estimates (\ref{mainintr1}), (\ref{mainintr2}) for both  Dirichlet and Neumann boundary conditions;
in Appendix we discuss some properties of the atlas distance $d_{{\mathcal{A}}}$,  the Hausdorff-Pompeiu lower deviation
$d_{{\mathcal{H}}{\mathcal{P}}}$ and the Hausdorff-Pompeiu distance
$d^{{\mathcal{H}}{\mathcal{P}}}$.

\section{Preliminaries and notation}
\label{section2}

Let $N, m \in \mathbb{N}$ and $\Omega$ be an open set in $\mathbb{R}^N$.
We denote by $W^{m,2}(\Omega )$ the Sobolev space  of
 complex-valued functions in $L^2(\Omega )$, which
have all distributional derivatives up to order $m$ in $L^2(\Omega ) $,
endowed with the norm
\begin{equation}
\label{norm}
\|u \|_{W^{m,2}(\Omega )}=\sum_{|\alpha|\le m}\|D^{\alpha }u \|_{L^{2}(\Omega )}.
\end{equation}
We denote by  $W^{m,2}_0(\Omega )$
the closure in $W^{m,2}(\Omega )$ of the space of the $C^{\infty}$-functions
with compact support in $\Omega$.\\

\begin{lem}
\label{sobolev}
Let $\Omega $ be an open set in ${\mathbb{R}}^N$. Let $V(\Omega )$ be a subspace of $W^{m,2}(\Omega )$ such that  the embedding $V(\Omega )\subset W^{m-1,2}(\Omega )$ is compact.
Then there exists $c>0$ such that
\begin{equation}
\label{equivalent}
\|u \|_{W^{m,2}(\Omega )}  \le c \biggl( \| u \|_{L^2(\Omega )} +\sum_{|\alpha|= m}\|D^{\alpha }u \|_{L^{2}(\Omega )}  \biggr),
\end{equation}
for all $u\in V(\Omega )$.
\end{lem}

{\bf Proof.} Since $(V(\Omega ),\| \cdot \|_{m,2} )$ is compactly embedded in $ W^{m-1,2}(\Omega )$ and $ W^{m-1,2}(\Omega )$ is continuously embedded in $L^2(\Omega )$,
by Lions' Lemma (cf. {\it e.g.}, Berger~\cite[p.~35]{ber}) it follows that for all $\epsilon \in ]0,1[$ there exists $c(\epsilon )>0$ such that
\begin{equation}
\| u\|_{W^{m-1,2}(\Omega )}\le \epsilon \| u\|_{W^{m,2}(\Omega )}+c(\epsilon )\| u\|_{L^{2}(\Omega )} \nonumber
\end{equation}
hence
\begin{equation}
\| u\|_{W^{m-1,2}(\Omega )}\le \frac{\epsilon}{(1-\epsilon )} \sum_{|\alpha|= m}\|D^{\alpha }u \|_{L^{2}(\Omega )}+\frac{c(\epsilon )}{(1-\epsilon )} \| u\|_{L^{2}(\Omega )},
\end{equation}
for all $u\in V(\Omega )$. Inequality (\ref{equivalent}) immediately follows.
\hfill $\Box$\\

Let $\hat m$ be the number of the
multi-indices $\alpha=(\alpha_1,\dots ,\alpha_N)\in {\mathbb{N}}_0^N$
with length $|\alpha | = \alpha_1+\dots +\alpha_N$ equal to $m$. Here ${\mathbb{N}}_0={\mathbb{N}}\cup\{0\} $.
For all $\alpha ,\beta \in
{\mathbb{N}}_0^N$
such that $|\alpha |=|\beta |=m$, let
$A_{\alpha \beta }$ be bounded measurable real-valued
functions defined on $\Omega$ satisfying
$A_{\alpha \beta }=
A_{\beta \alpha }$ and the uniform ellipticity
condition
\begin{equation}
\label{elp}
\sum_{|\alpha |=|\beta |=m}A_{\alpha \beta }(x)\xi_{\alpha }\xi_{\beta }
\geq \theta |\xi |^2
\end{equation}
for all $x\in \Omega $, $\xi =(\xi_{\alpha })_{|\alpha |=m}\in\mathbb{R}^{\hat m}$,
where $\theta >0$ is the ellipticity constant.

Let $V(\Omega )$ be a closed subspace of  $W^{m,2}(\Omega )$ containing  $W^{m,2}_0(\Omega )$. We consider the
following eigenvalue problem

\begin{equation}
\label{mainprobl}
\int_{\Omega}\sum_{|\alpha |=|\beta |=m}
A_{\alpha \beta }D^{\alpha }u D^{\beta }\bar v dx=\lambda \int_{\Omega}u \bar v dx,
\end{equation}
for all test functions $v\in V (\Omega )$, in the unknowns
$u\in V (\Omega )$ (the eigenfunctions) and $\lambda \in\mathbb{R}$
(the eigenvalues).

Clearly problem (\ref{mainprobl}) is the weak formulation of an eigenvalue problem for the operator
$H$ in (\ref{classic}) subject to suitable homogeneous boundary conditions and the choice of $V(\Omega )$ corresponds to the choice of the boundary conditions (see {\it e.g.,} Ne\v{c}as \cite{Ne67}).

We set
\begin{equation}
Q_{\Omega }(u,v)=
\int_{\Omega }
\sum_{|\alpha |=|\beta |=m} A_{\alpha \beta }D^{\alpha }u
D^{\beta }\bar vdx,\ \ \ \ \ \ \ Q_{\Omega }(u)=Q_{\Omega }(u,u)  ,
\end{equation}
 for all $u,v\in W^{m,2}(\Omega ) $.

If the embedding $V(\Omega )\subset W^{m-1,2}(\Omega )$ is compact, then
the eigenvalues of equation (\ref{mainprobl}) coincide with the eigenvalues of
a suitable operator $H_{V(\Omega ) }$ canonically associated with the restriction of the
quadratic form $Q_{\Omega }$ to $V(\Omega )$.
In fact, we have the following theorem.

\begin{thm}
\label{setting}
Let $\Omega$ be an open set in $\mathbb{R}^N$.
Let $m\in\mathbb{N}$, $\theta >0$ and, for all $\alpha ,\beta \in
 {\mathbb{N}}_0^{N} $
such that $|\alpha |=|\beta |=m$, let
$A_{\alpha \beta }$ be bounded measurable real-valued
functions defined on $\Omega$,  satisfying
$A_{\alpha \beta }= A_{ \beta\alpha  }$ and condition (\ref{elp}).

Let $V (\Omega )$ be a closed subspace of  $W^{m,2}(\Omega )$ containing  $W^{m,2}_0(\Omega )$  and such that the embedding
$V(\Omega )\subset W^{m-1,2}(\Omega )$ is compact.

Then there exists a non-negative self-adjoint linear operator
$H_{V(\Omega )}$ on $L^2(\Omega )$ with compact resolvent,
 such that ${\rm Dom}(H^{1/2}_{V(\Omega )})=
V (\Omega )$ and
\begin{equation}
\label{setting0}
<H^{1/2}_{V(\Omega )}u,H_{V(\Omega )}^{1/2} v >_{L^2(\Omega )}=
Q_{\Omega }(u,v)
,
\end{equation}
for all  $u,v \in V(\Omega )$.
Moreover, the eigenvalues of equation
(\ref{mainprobl})
coincide with the eigenvalues
$\lambda_n [H_{V(\Omega )} ]$ of
$H_{V(\Omega )} $ and
\begin{equation}
\label{setting1}
\lambda_n [H_{V(\Omega )}]
=\inf_{\substack{{\mathcal{L}}\le   V (\Omega ) \\ {\rm dim}\, {\mathcal{L}}=n}}\sup_{\substack{u\in {\mathcal{L}} \\ u\ne 0}}
\frac{Q_{\Omega }(u)}{\|u \|_{L^2(\Omega )}^2}
 .
\end{equation}
\end{thm}

{\bf Proof.}
By Lemma \ref{sobolev}, inequality (\ref{elp}) and by the
boundedness of the coefficients $A_{\alpha \beta}$, it follows that
the space $V(\Omega )$ endowed with the norm defined by
\begin{equation}
\label{setting4}
(\| u\|^2_{L^2(\Omega )}+Q_{\Omega }(u))^{1/2},
\end{equation}
for all $u\in V(\Omega )$, is complete. Indeed, this norm is equivalent on $V(\Omega )$ to the norm defined by (\ref{norm}).
Thus, the restriction of the quadratic form
$Q_{\Omega }$ to $V(\Omega )$ is a closed quadratic form on $V(\Omega )$ (cf. Davies~\cite[pp.~81-83]{da}) and there exists a non-negative
self-adjoint operator $H_{V(\Omega )}$ on $L^2(\Omega )$ satisfying ${\rm Dom}(H_{V(\Omega )}^{1/2})=
V(\Omega )$ and condition (\ref{setting0}) (cf. \cite[Theorem~4.4.2]{da}).
Since  the embedding
$V(\Omega )\subset L^2(\Omega )$ is compact then $H_{V(\Omega )}$ has
compact resolvent (cf. \cite[Ex.~4.2]{da}). The fact that the eigenvalues of equation (\ref{mainprobl}) coincide with the eigenvalues
of the operator $H_{V(\Omega )} $
is well known. Finally, the variational representation in (\ref{setting1}) is
given by the well-known Min-Max Principle (cf. \cite[Theorem~4.5.3]{da}). \hfill $\Box$

\begin{defn}\label{defset} Let $\Omega$ be an open set in ${\mathbb{R}}^N$.
Let $m\in\mathbb{N}$, $\theta >0$ and, for all $\alpha ,\beta \in
 {\mathbb{N}}_0^{N} $
such that $|\alpha |=|\beta |=m$, let
$A_{\alpha \beta }$ be bounded measurable real-valued
functions defined on $\Omega$,  satisfying
$A_{\alpha \beta }= A_{ \beta\alpha  }$ and condition (\ref{elp}).

If the embedding $W^{m,2}_0(\Omega )\subset  W^{m-1,2}(\Omega )$ is compact, we set
$$
\lambda_{n,{\mathcal {D}}}[\Omega ]=\lambda_n[H_{W^{m,2}_0(\Omega )}  ].
$$

If the embedding $W^{m,2}(\Omega )\subset  W^{m-1,2}(\Omega )$ is compact,
we set
$$
\lambda_{n,{\mathcal {N}}}[\Omega ]=\lambda_n[H_{W^{m,2}(\Omega )}  ].
$$

The numbers $\lambda_{n,{\mathcal {D}}}[\Omega ]$, $\lambda_{n,{\mathcal {N}}}[\Omega ]$   are called the Dirichlet eigenvalues, Neumann eigenvalues
respectively,  of the operator (\ref{classic}).
\end{defn}

When we refer to both Dirichlet and Neumann boundary conditions we write just $\lambda_n[\Omega ]$ instead of
$\lambda_{n,{\mathcal {D}}}[\Omega ]$ and $\lambda_{n,{\mathcal {N}}}[\Omega ]$.

\begin{rem}
\label{sufcomp}
If $\Omega $ is such that the embedding $W^{1,2}_0(\Omega )\subset L^2(\Omega )$ is compact (for instance, if $\Omega $ is an arbitrary open set with  finite Lebesgue measure) then also the embedding $ W^{m,2}_0(\Omega ) $ $  \subset W^{m-1,2}(\Omega )$ is compact and the Dirichlet eigenvalues are well-defined.

If $\Omega $ is such that the embedding $W^{1,2}(\Omega )\subset L^2(\Omega )$ is compact (for instance, if $\Omega $
has a continuous boundary, see Definition~\ref{class}) then the embedding $W^{m,2}(\Omega )\subset W^{m-1,2}(\Omega )$ is compact and
the Neumann eigenvalues are well-defined.
\end{rem}

\begin{example}

Let $\Omega$ be an open set in $ {\mathbb{R}}^2 $. We consider the bi-harmonic operator $\Delta ^2$ in ${\mathbb{R}}^2$ and the sesquilinear form
$$
Q_{\Omega }(u,v)=\int_{\Omega }\biggl(\frac{\partial ^2u}{\partial x_1^2}\frac{\partial ^2\bar v}{\partial x_1^2}+2\frac{\partial ^2 u}{\partial x_1\partial x_2} \frac{\partial ^2 \bar v}{\partial x_1\partial x_2}+\frac{\partial ^2u}{\partial x_2^2}\frac{\partial ^2\bar v}{\partial x_2^2}\biggr)dx,\ \ u,v\in V(\Omega ),
$$
where $V(\Omega )$ is either $W^{2,2}_0(\Omega)$ (Dirichlet boundary conditions) or $W^{2,2}(\Omega)$ (Neumann boundary conditions). Recall that
the Euler-Lagrange equation for the minimization of the quadratic form $Q_{\Omega }(u,u)$ is $\Delta^2u=0$.
Observe that condition (\ref{elp}) is satisfied with $\theta =1$.

Let $H_{V(\Omega )}$ be the operator associated with the quadratic form $Q_{\Omega }$ as in Theorem~\ref{setting}. Consider the eigenvalue problem
\begin{equation}
\label{biha}
H_{V(\Omega )}u=\lambda u .
\end{equation}

In the case $V(\Omega )=W^{2,2}_0(\Omega )$  equation (\ref{biha}) is the weak formulation of the classical eigenvalue problem for the bi-harmonic operator
subject to Dirichlet  boundary conditions
\begin{equation}
\label{bihadir}\left\{
\begin{array}{ll}
\Delta ^2u=\lambda u, & {\rm in }\ \Omega ,\vspace{1mm} \\
u=0, &  {\rm on }\ \partial\Omega ,\vspace{1mm}\\
\frac{\partial u}{\partial n}=0, &  {\rm on }\ \partial\Omega ,
\end{array}\right.
\end{equation}
for bounded domains $\Omega$ of class $C^2$. Here $n =(n_1, n_2)$ is the unit outer normal to $\partial \Omega $.

In the case $V(\Omega )=W^{2,2}(\Omega )$  equation (\ref{biha}) is the weak formulation of the classical eigenvalue problem for the bi-harmonic operator
subject to Neumann  boundary conditions

\begin{equation}\label{bihaneu}\left\{
\begin{array}{ll}
\Delta ^2u=\lambda u, & {\rm in }\ \Omega ,\vspace{1mm} \\
\frac{\partial^2 u}{\partial n ^2}=0, &  {\rm on }\ \partial\Omega ,\vspace{1mm}\\
\frac{d}{ds}\frac{\partial^2 u}{\partial n \partial t}+\frac{\partial \Delta u}{\partial n}=0, &  {\rm on }\ \partial\Omega ,
\end{array}\right.
\end{equation}
for bounded domains $\Omega$ of class $C^2$. Here $$\frac{\partial^2 u}{\partial n ^2}=\sum_{i,j=1}^2\frac{\partial ^2u}{\partial x_i\partial x_j}n_i n_j,\ \ \ \frac{\partial^2 u}{\partial n \partial t}=\sum_{i,j=1}^2\frac{\partial ^2u}{\partial x_i\partial x_j}n_i t_j,$$
 $s$ denotes the arclengh  of $\partial \Omega  $ (with positive orientation), $t=(t_1,t_2)$ denotes the unit tangent vector  to $\partial \Omega $ (oriented in the sense of increasing $s$). This follows by a standard argument
and by observing that if   $u,v\in C^4(\bar \Omega )$  then by  the Divergence Theorem
\begin{eqnarray*}
Q_{\Omega }(u,v)= \int_{\Omega }\Delta^2 u\bar v dx-\int_{\partial \Omega }\frac{\partial \Delta u}{\partial n}\bar v d\sigma +
\int_{\partial \Omega }\left(n_{1}\nabla \frac{\partial u}{\partial x_1}+n_{2}\nabla \frac{\partial u}{\partial x_2}\right)\cdot \nabla \bar vd\sigma \\
= \int_{\Omega }\Delta^2 u \bar vdx+\int_{\partial \Omega }\left(\frac{\partial^2 u}{\partial n ^2}\frac{\partial\bar v}{\partial n}+\frac{\partial^2u}{\partial n \partial t }\frac{\partial\bar v}{\partial t}-\frac{\partial \Delta u}{\partial n}\bar v \right) d\sigma .
\end{eqnarray*}

One may also consider the sesquilinear form
$$
Q_{\Omega }^{(\nu )}(u,v )=\nu  \int_{\Omega }\Delta u\Delta \bar v  + (1-\nu )Q_{\Omega }(u,v),\ \ u, v\in V(\Omega ).
$$
If $0\le \nu <1$ then condition (\ref{elp}) is satisfied with $\theta =1-\nu$. Observe that
the Euler-Lagrange equation for the minimization of the quadratic form $Q_{\Omega }^{(\nu )}(u,u)$ is again $\Delta^2u=0$.

Let $H_{V(\Omega )}^{(\nu )}$ be the operator associated with the quadratic form $Q_{\Omega }^{(\nu )}$ as in Theorem~\ref{setting}. Consider the eigenvalue problem
\begin{equation}
\label{bihabis}
H_{V(\Omega )}^{(\nu )}u=\lambda u .
\end{equation}

In the case $V(\Omega )=W^{2,2}_0(\Omega )$  equation (\ref{bihabis}) is another weak formulation of the classical eigenvalue problem (\ref{bihadir}).

In the case $V(\Omega )=W^{2,2}(\Omega )$  equation (\ref{bihabis}) is the weak formulation
of the classical eigenvalue problem for the bi-harmonic operator
subject to Neumann  boundary conditions depending on $\nu $
\begin{equation}\label{bihaneubis}\left\{
\begin{array}{ll}
\Delta ^2u=\lambda u, & {\rm in }\ \Omega ,\vspace{1mm} \\
\nu \Delta u +(1-\nu )\frac{\partial^2 u}{\partial n ^2}=0, &  {\rm on }\ \partial\Omega ,\vspace{1mm}\\
(1-\nu )\frac{d}{ds}\frac{\partial^2 u}{\partial n \partial t}+\frac{\partial \Delta u}{\partial n}=0, &  {\rm on }\ \partial\Omega .
\end{array}\right.
\end{equation}

The bi-harmonic operator subject to these boundary conditions with $0<\nu <1/2$ arises in the study of
 small deformations of a thin plate
under Kirchhoff  hypothesis in which case $\nu$ is the Poisson ratio of the plate (see {\it e.g.}, Nazaret~\cite{naza} and the references therein).

\end{example}

\section{Open sets with continuous boundaries}
\label{section3}

We recall that for any set $V$ in ${\mathbb{R}}^N$ and $\delta >0$ we denote by $V_{\delta }$ the set $\{x\in V:\ d(x, \partial \Omega )>\delta \}$. Moreover, by a rotation in ${\mathbb{R}}^N$ we mean a $N\times N$-orthogonal matrix with real entries which we identify with the corresponding linear operator acting in ${\mathbb{R}}^N$.

\begin{defn}
\label{class}

Let $ \rho >0$, $s,s'\in\mathbb{N}$, $s'\le s$
and  $\{V_j\}_{j=1}^s$ be a family of bounded open cuboids  and
$\{r_j\}_{j=1}^{s} $ be a family of rotations in ${\mathbb{R}}^N $.

We say that that ${\mathcal{A}}= (  \rho , s,s', \{V_j\}_{j=1}^s, \{r_j\}_{j=1}^{s} ) $ is an atlas in ${\mathbb{R}}^N$ with the parameters
$\rho , s,s', \{V_j\}_{j=1}^s, \{r_j\}_{j=1}^{s}$, briefly an atlas in ${\mathbb{R}}^N$.

We denote by $C( {\mathcal{A}}   )$ the family of all open sets $\Omega $ in ${\mathbb{R}}^N$
satisfying the following properties:

(i) $ \Omega\subset \bigcup\limits_{j=1}^s(V_j)_{\rho}$ and $(V_j)_\rho\cap\Omega\ne\emptyset;$

(ii) $V_j\cap\partial \Omega\ne\emptyset$ for $j=1,\dots s'$, $ V_j\cap \partial\Omega =\emptyset$ for $s'<j\le s$;

(iii) for $j=1,...,s$
$$
r_j(V_j)=\{\,x\in \mathbb{R}^N:~a_{ij}<x_i<b_{ij}, \,i=1,....,N\},
$$

\noindent and

$$
r_j(\Omega\cap V_j)=\{x\in\mathbb{R}^N:~a_{Nj}<x_N<g_{j}(\bar
x),~\bar x\in W_j\},$$

\noindent where $\bar x=(x_1,...,x_{N-1})$, $W_j=\{\bar
x\in\mathbb{R}^{N-1}:~a_{ij}<x_i<b_{ij},\,i=1,...,N-1\}$
and $g_j$ is a continuous function defined on $\overline {W}_j$ (it is meant that if $s'<j\le s$ then $g_j(\bar x)=b_{Nj}$ for all $\bar x\in \overline{W}_j$);

moreover for $j=1,\dots ,s'$
$$
a_{Nj}+\rho\le g_j(\bar x)\le b_{Nj}-\rho ,$$

\noindent for all $\bar x\in \overline{W}_j$.

We say that an open set $\Omega$ in ${\mathbb{R}}^N$ is an open set with a continuous boundary if $\Omega $ is of class
$C( {\mathcal{A}}   ) $ for some atlas $ {\mathcal{A}}   $.

\end{defn}

We note that, for an open set $\Omega$ of class $C({\mathcal{A}}) $, inequality (\ref{equivalent}) holds for all $u\in W^{m,2}(\Omega ) $ with
a constant $c$ depending only on ${\mathcal{A}}$. More precisely, we denote by ${\mathcal{D}}_{\Omega } $  {\it the  best constant} for which inequality (\ref{equivalent}) is satisfied
for $V(\Omega )=W^{m,2}_0(\Omega ) $.
We denote by  ${\mathcal{N}}_{\Omega } $  {\it the best constant} for which inequality (\ref{equivalent}) is satisfied
for $V(\Omega )=W^{m,2}(\Omega ) $. Then we have the following (for a proof we refer to Burenkov~\cite[Thm. 6, p.~160]{bu}).

\begin{lem}
\label{interm}
Let ${\mathcal{A}} $ be an atlas in ${\mathbb{R}}^N$, $m\in {\mathbb{N}}$. There exists $c>0$ depending only on $N, {\mathcal{A} }$ and $m$ such that
\begin{equation}
1\le {\mathcal{D}}_{\Omega }\le {\mathcal{N}}_{\Omega }\le c ,
\end{equation}
for all open sets $\Omega\in C({\mathcal{A}}) $.
\end{lem}

\begin{lem}
\label{unifbound}
Let ${\mathcal{A}}$ be an atlas in ${\mathbb{R}}^N$. Let $  m\in {\mathbb{N}}$,  $L, \theta >0$ and, for all $\alpha ,\beta \in {\mathbb{N}}_0^N$ with $|\alpha |=|\beta |=m$, let $A_{\alpha \beta }\in L^{\infty }(\cup_{j=1}^sV_j)$  satisfy $A_{\alpha \beta }=A_{ \beta \alpha }$, $\| A_{\alpha\beta} \|_{ L^{\infty }(\cup_{j=1}^sV_j)}\le L$ and condition (\ref{elp}).

Then for each $n\in {\mathbb{N}}$ there exists $\Lambda_n>0$ depending only on $n, N, {\mathcal{A}}, m$ and  $L$  such that
\begin{equation}
\label{newnumb}
\lambda_{n,{\mathcal{N}}}[\Omega ]\le \lambda_{n,{\mathcal{D}}}[\Omega ]\le \Lambda_n ,
\end{equation}
for all open sets $\Omega \in C({\mathcal{A}}) $.
\end{lem}

{\bf Proof.} The inequality $\lambda_{n,{\mathcal{N}}}[\Omega ]\le \lambda_{n,{\mathcal{D}}}[\Omega ]$  is well known. Now we prove the second inequality. Clearly, there exists a ball $B$ of radius $\rho /2$ such that $B\subset \Omega $ for all open sets $\Omega \in C({\mathcal{A}}) $.
By the well-known monotonicity of the Dirichlet eigenvalues with respect to inclusion it follows that
$$
\lambda_{n,{\mathcal{D}}}[\Omega ]\le \lambda_{n,{\mathcal{D}}}[B ].
$$
Thus it suffices to estimate  $\lambda_{n,{\mathcal{D}}}[B ]$.
Clearly there exists $c>0$ depending only on $N$ and $m $ such that
\begin{equation}
\label{monun}
Q_{B}(u)\le cL\int_{B}|\nabla ^m u|^2dx,
\end{equation}
for all $u\in W^{m,2}_0(B)$, where $\nabla ^m u=(D^{\alpha }u)_{|\alpha |=m}  $. By (\ref{setting1}) and (\ref{monun}) it follows that
$$
\lambda_{n,{\mathcal{D}}}[B]\le\Lambda_n\equiv  cL \inf_{\substack{{\mathcal{L}}\le W^{m,2}_0(B)\\ {\rm dim}\, {\mathcal{L}}=n  }} \sup_{\substack{u\in {\mathcal{L}}\\ u\ne 0}} \frac{\int_{B}|\nabla ^m u|^2dx}{\int_B|u|^2dx}<\infty .
$$
Clearly $\Lambda_n$ depends only on $n, N,\rho,  m $ and $L$.
 \hfill $\Box$

\section{The case of diffeomorphic open sets}
\label{section4}

\begin{lem}
\label{variables}
 Let $\Omega$ be an  open set in $\mathbb{R}^N$. Let $m\in {\mathbb{N}}$, $B_1,B_2>0$ and $\phi$ be a diffeomorphism of $\Omega$ onto $\phi(\Omega )$ of class
$C^m$ such that
\begin{equation}
\label{varia1}
\max_{1\le |\alpha |\le m} |D^{\alpha }\phi (x)
| \le B_1, \ \ \ \  |{\rm det }\nabla \phi (x)|\geq  B_2,
\end{equation}
for all $x\in \Omega $.
Let $B_3>0$ and, for all $\alpha ,\beta \in
\mathbb{N}_0^N$
such that $|\alpha |=|\beta |=m$, let
$A_{\alpha \beta }$ be measurable real-valued
functions defined on $\Omega\cup\phi (\Omega )$ satisfying
\begin{equation}
\label{varia2}
\max_{|\alpha |= |\beta |= m}|A_{\alpha \beta }(x) |\le B_3,
\end{equation}
for almost all $x\in \Omega \cup \phi (\Omega )$.
Then there exists $c>0$ depending only on $N,
m, B_1,B_2,B_3$ such that

\begin{equation}
\label{num1}
\left| Q_{\phi (\Omega )}(u\circ\phi ^{(-1)})-
Q_{\Omega }(u) \right| \le c{\mathcal{L}}(\phi) \int_{\Omega}
\sum_{1\le |\alpha|\le m}|D^{\alpha }u|^2dx,
\end{equation}
for all $u\in W^{m,2}(\Omega )$, where
\begin{equation}
\label{variables2}
{\mathcal{L}}(\phi )= \max_{1\le |\alpha |\le m}   \| D^{\alpha }(\phi -{\rm Id} )\|_{L^{\infty }(\Omega )}
 +\max_{|\alpha |=|\beta |=m}
\| A_{\alpha \beta }\circ \phi -A_{\alpha \beta } \|_{L^{\infty}(\Omega )}.
\end{equation}
\end{lem}

{\bf Proof.}
By changing variables and using a known formula for high derivatives
of composite functions (cf. {\it e.g.} Fraenkel \cite[Formula~B]{fr}),
we have that
\begin{eqnarray}
 \label{frae}
\lefteqn{
Q_{\phi (\Omega)}(u\circ \phi^{(-1)})=
\int_{\phi (\Omega )}\sum_{|\alpha |=|\beta |=m}A_{\alpha \beta}
D^{\alpha }(u\circ \phi^{(-1)})\, \overline{ D^{\beta }(u\circ \phi^{(-1)})}dy
}  \nonumber \\
& & \quad=
\int_{\Omega }\sum_{|\alpha |=|\beta |=m} \left(A_{\alpha \beta}
D^{\alpha }(u\circ \phi^{(-1)})
\overline{D^{\beta }(u\circ \phi^{(-1)}   )}\right)\circ\phi
|{\rm det }\nabla \phi |dx \nonumber \\
& &  \quad =
\int_{\Omega }\sum_{|\alpha |=|\beta |=m} A_{\alpha \beta}\circ \phi
\sum_{\substack{ 1\le |\eta |\le |\alpha | \\
1\le |\xi |\le |\beta |  } } D^{\eta }u\overline{D^{\xi }u}\,\,
(
{p}_{\alpha\eta}(\phi^{(-1)})
{p}_{\beta\xi}(\phi^{(-1)}))\circ\phi
|{\rm det }\nabla \phi |dx
\nonumber   \\ & &   \quad =
\sum_{\substack{|\alpha |=|\beta |=m\\
 1\le |\eta |\le |\alpha | \\
 1\le |\xi |\le |\beta |
}}
\int_{\Omega }(A_{\alpha \beta}
{p}_{\alpha\eta}(\phi^{(-1)})
{p}_{\beta\xi}(\phi^{(-1)}))\circ\phi\,
D^{\eta }u\overline{D^{\xi }u}\, |{\rm det }\nabla \phi|  dx ,
\end{eqnarray}
for all $u\in W^{m,2} (\Omega)$, where
for all $\alpha , \eta$ with  $1\le |\eta |\le |\alpha|= m$,
${p}_{\alpha\eta}(\phi^{(-1)}) $
denotes a polynomial  of degree $|\eta |$ in derivatives
of $\phi^{(-1)}$ of order between $1$ and $|\alpha |$, with coefficients depending only on
$N,\alpha ,\eta $.

We recall that for each $\alpha$ with  $1\le |\alpha|\le m$ there exists a polynomial
$p_{\alpha}(\phi )$
in derivatives of $\phi$ of order between $1$ and $|\alpha |$,  with coefficients depending only on $N, \alpha$, such that
\begin{equation}
\label{frae1}
(D^{\alpha }\phi^{(-1)}) \circ \phi =\frac{p_{\alpha}(\phi )}{({\rm det}\nabla \phi )^{2{|\alpha |}-1}}.
\end{equation}

In order to estimate $Q_{\phi (\Omega )}(u\circ\phi ^{(-1)})-
Q_{\Omega }(u)$ it is enough to estimate the expressions
\begin{equation}
(A_{\alpha \beta}
{p}_{\alpha\eta}(\phi^{(-1)})
{p}_{\beta\xi}(\phi^{(-1)}))\circ\phi\,
\, |{\rm det }\nabla \phi| -(A_{\alpha \beta}
{p}_{\alpha\eta}(\tilde \phi^{(-1)})
{p}_{\beta\xi}(\tilde \phi^{(-1)}))\circ\tilde \phi\,
\, |{\rm det }\nabla \tilde \phi| ,\nonumber
\end{equation}
where $\tilde\phi ={\rm Id}$. This can be done by using the triangle inequality and by observing that (\ref{frae1}) implies that
$$
|(D^{\alpha }\phi^{(-1)}) \circ \phi -(D^{\alpha }\tilde \phi^{(-1)}) \circ\tilde  \phi |\le c
 \max_{1\le |\beta |\le |\alpha |}   \| D^{\beta }(\phi -\tilde \phi )\|_{L^{\infty }(\Omega )},
$$
where $c$ depends only on $N, \alpha , B_1, B_2$.
\hfill $\Box $\\

\begin{thm}
\label{diffeothm}Let $U$ be an open set in ${\mathbb{R}}^N$.
 Let $m\in {\mathbb{N}}$, $B_1,B_2, B_3, \theta >0$.
For all $\alpha ,\beta \in
 \mathbb{N}_0^{N} $
with $|\alpha |=|\beta |=m$, let
$A_{\alpha \beta }$ be measurable real-valued
functions defined on $U$,  satisfying
$A_{\alpha \beta }= A_{ \beta\alpha  }$ and conditions (\ref{elp}), (\ref{varia2}) in $U$. The following statements hold.
\begin{itemize}
\item[(i)]
There exists $c_1>0$ depending only on $N,m , B_1 ,B_2, B_3, \theta $ such that for all $n\in {\mathbb{N}}$,
for all  open sets $\Omega\subset U$ such that the embedding $W^{m,2}_0(\Omega)\subset W^{m-1,2}(\Omega)$ is compact,  and for all diffeomorphisms of $\Omega$ onto $\phi(\Omega )$ of class
$C^m$ satisfying (\ref{varia1}) and such that $\phi (\Omega )\subset U$, the inequality
\begin{equation}
\label{diffeothm1}
|\lambda_{n, {\mathcal{D}}}[\Omega ]-\lambda_{n, {\mathcal{D}}}[\phi (\Omega )]|\le c_1\, {\mathcal{D}}^2_{\Omega }(1+\lambda_{n, {\mathcal{D}}}[\Omega ]) {\mathcal{L}}(\phi )
\end{equation}
holds  if $ {\mathcal{L}}(\phi )< (c_1\, {\mathcal{D}}^2_{ \Omega })^{-1} $.
\item[(ii)]
There exists $c_2>0$ depending only on $N,m, B_1 ,B_2, B_3, \theta $ such that for all $n\in {\mathbb{N}}$,
for all  open sets $\Omega\subset U$ such that the embedding $W^{m,2}(\Omega)\subset W^{m-1,2}(\Omega)$ is compact, and for all diffeomorphisms of $\Omega$ onto $\phi(\Omega )$ of class
$C^m$ satisfying (\ref{varia1}) and such that $\phi (\Omega )\subset U$, inequality

\begin{equation}
\label{diffeothm2}
|\lambda_{n, {\mathcal{N}}}[\Omega ]-\lambda_{n, {\mathcal{N}}}[\phi (\Omega )]|\le c_2\, {\mathcal{N}}^2_{\Omega }(1+\lambda_{n, {\mathcal{N}}}[\Omega ]) {\mathcal{L}}(\phi ),
\end{equation}
holds if $ {\mathcal{L}}(\phi )< (c_2\, {\mathcal{N}}^2_{ \Omega })^{-1} $.
\end{itemize}
\end{thm}

{\bf Proof.} We prove statement $(i)$. Let $\Omega\subset U$ be an open set
such that the embedding $W^{m,2}_0(\Omega)\subset W^{m-1,2}(\Omega)$ is compact
and $\phi$ be a diffeomorphisms of $\Omega$ onto $\phi(\Omega )$ of class
$C^m$ satisfying (\ref{varia1}) and such that $\phi (\Omega )\subset U$.
By inequalities (\ref{equivalent}), (\ref{elp}), (\ref{num1})   it follows that there exists $c_3>0$ depending only on $N,m, B_1 ,B_2, B_3, \theta $
such that
\begin{equation}
\label{num1'}
\left| Q_{\phi (\Omega )}(u\circ\phi ^{(-1)})-
Q_{\Omega }(u) \right| \le c_3 \, {\mathcal{D}}^2_{\Omega } ( \| u \|^2_{L^2(\Omega )} +Q_{\Omega }(u) ) {\mathcal{L}}(\phi).
\end{equation}

Clearly we have

\begin{eqnarray}
\label{tria}
\lefteqn{\left|
\frac{Q_{\phi (\Omega )}(u\circ\phi ^{(-1)})}{\|u\circ\phi ^{(-1)} \|_
{L^2(\phi (\Omega ))}^2
}
-
\frac{Q_{\Omega }(u)}{\|u \|_
{L^2(\Omega )}^2
}
\right|} \nonumber \\
& & \le
\frac{| Q_{\phi (\Omega )}(u\circ\phi ^{(-1)})- Q_{\Omega }(u) |
}{\int_{\Omega}|u|^2|{\rm det}\nabla \phi | dx } +
\frac{Q_{\Omega }(u)\int_{\Omega}|u|^2 ||{\rm det}\nabla \phi |-1| dx  }
{\int_{\Omega}|u|^2|{\rm det}\nabla \phi | dx\int_{\Omega }|u|^2dx }.
\end{eqnarray}
By observing that ${\mathcal{D}}_{\Omega }\geq 1$ and by combining inequalities (\ref{num1'}) and  (\ref{tria}) it follows that there exists
$c_4>0$ depending only on  $N,m, B_1 ,B_2, B_3, \theta $
such that for all $u\in W^{m,2}_0(\Omega )$
\begin{equation}
\left|
\frac{Q_{\phi (\Omega )}(u\circ\phi ^{(-1)})}{\|u\circ\phi ^{(-1)} \|_
{L^2(\phi (\Omega ))}^2
}
-
\frac{Q_{\Omega }(u)}{\|u \|_
{L^2(\Omega )}^2
}
\right|\le c_4\, {\mathcal{D}}^2_{\Omega }\left(1+ \frac{Q_{\Omega }(u)}{\|u \|_
{L^2(\Omega )}^2
} \right){\mathcal{L}}(\phi ),
\end{equation}
which can be written as
\begin{eqnarray}
\label{minmaxappl}
\lefteqn{(1-c_4{\mathcal{D}}^2_{\Omega }{\mathcal{L}}(\phi ) )\frac{Q_{\Omega }(u)}{\|u \|_
{L^2(\Omega )}^2}-c_4{\mathcal{D}}^2_{\Omega } {\mathcal{L}}(\phi )} \\  \nonumber & & \qquad\qquad\qquad\qquad\qquad \le
\frac{Q_{\phi (\Omega )}(u\circ\phi ^{(-1)})}{\|u\circ\phi ^{(-1)} \|_
{L^2(\phi (\Omega ))}^2
}\\ \nonumber & & \qquad \qquad\qquad\qquad\qquad \le (1+c_4{\mathcal{D}}^2_{\Omega }{\mathcal{L}}(\phi ) )\frac{Q_{\Omega }(u)}{\|u \|_
{L^2(\Omega )}^2}+c_4{\mathcal{D}}^2_{\Omega } {\mathcal{L}}(\phi ).
\end{eqnarray}
Assume now that  $1-c_4{\mathcal{D}}^2_{\Omega }{\mathcal{L}}(\phi )>0$.
Observe that the map $C_{\phi }$ of $L^2(\Omega )$ to  $L^2(\phi (\Omega ))$ which takes $u\in L^2(\Omega )$ to $C_{\phi }u=u\circ \phi^{-1}$ is a linear homeomorphism which restricts to a linear homeomorphism of $W^{m,2}_0(\Omega )$ onto $W^{m,2}_0(\phi (\Omega ) )$,
and that the embedding $ W^{m,2}_0(\phi (\Omega ))\subset W^{m-1,2}(\phi (\Omega )) $ is compact. Then
by applying the Min-Max Principle (\ref{setting1}) and using inequality  (\ref{minmaxappl}), it easy to deduce the validity of  inequality (\ref{diffeothm1}).

The proof of statement $(ii)$ is very  similar. In this case one should observe that the map $C_{\phi }$ defined above restricts to a linear homeomorphism
of $W^{m,2}(\Omega )$ onto $W^{m,2}(\phi (\Omega ) )$
and that if the embedding $W^{m,2}(\Omega )\subset W^{m-1,2}(\Omega ) $ is compact then also the embedding
$W^{m,2}(\phi (\Omega ))\subset W^{m-1,2}(\phi (\Omega ) ) $ is compact. \hfill $\Box $\\

\begin{corol}
\label{diffeocorol}
Let ${\mathcal{A}}$ be an atlas in ${\mathbb{R}}^N$. Let $m\in {\mathbb{N}}$,  $B_1, B_2, L,\theta >0$ and, for all $\alpha ,\beta \in {\mathbb{N}}_0^N$ with $|\alpha |=|\beta |=m$, let $A_{\alpha \beta }\in C^{0,1}(\cup_{j=1}^sV_j)$  satisfy $A_{\alpha \beta }=A_{ \beta \alpha }$, $\| A_{\alpha\beta} \|_{ C^{0,1}(\cup_{j=1}^sV_j)}\le L $ and condition (\ref{elp}).

Then there exists $c>0$ depending only on $N, {\mathcal{A}}, m, B_1 ,B_2, L, \theta  $ such that for all $n\in {\mathbb{N}}$,
for all  open sets $\Omega \in C({\mathcal{A}})$,   and for all diffeomorphisms of $\Omega$ onto $\phi(\Omega )$ of class
$C^m$ satisfying (\ref{varia1}) and such that $\phi (\Omega )\subset \cup_{j=1}^sV_j$, the inequality
\begin{equation}
\label{diffeocorol1}
|\lambda_{n}[\Omega ]-\lambda_{n}[\phi (\Omega )]|\le c(1+\lambda_{n}[\Omega ]) \max_{0\le |\alpha |\le m}   \| D^{\alpha }(\phi -{\rm Id} )\|_{L^{\infty }(\Omega )}
\end{equation}
holds for both Dirichlet and Neumann  boundary conditions, if $$ \max_{0\le |\alpha |\le m}   \| D^{\alpha }(\phi -{\rm Id} )\|_{L^{\infty }(\Omega )} < c^{-1} . $$
\end{corol}

{\bf Proof.}  It suffices to apply Lemma \ref{interm} and Theorem \ref{diffeothm}. \hfill $\Box$\\

\section{Estimates for Dirichlet eigenvalues via the atlas distance}
\label{section5}

\begin{defn}
\label{dev}
Let ${\mathcal{A}} =(\rho , s, s', \{V_j\}_{j=1}^s , \{r_j\}_{j=1}^s )$ be an atlas in ${\mathbb{R}}^N$.
For all $\Omega_1 , \Omega_2\in C({\mathcal{A}})$ we define the `atlas distance' $d_{{\mathcal{A}}}$ by

\begin{equation}
d_{{\mathcal{A}}}(\Omega_1,\Omega_2)=\max_{j=1, \dots ,s}\sup_{(\bar x , x_N)\in r_j(V_j)}\left| g_{1j}(\bar x) - g_{2j}(\bar x)  \right|,
\end{equation}
where $g_{1j}$, $g_{2j}$ respectively, are the functions describing the boundaries of $\Omega_1, \Omega_2$ respectively, as in Definition \ref{class} $(iii)$.
\end{defn}

We observe that the function $d_{{\mathcal{A}}}(\cdot, \cdot)$ is in fact a distance in $C({{\mathcal{A}}} )$ (for further properties of $d_{{\mathcal{A}}}$ see also Appendix).

If
$\Omega\in C({\mathcal{A}})$ it will be useful to set
\begin{equation}
\label{dj}
d_j(x,\partial \Omega
)= |g_j(\, \overline{(r_j(x))}\, )-(r_j(x))_N |,
\end{equation}
for all  $j=1, \dots , s$ and $x\in V_j$, where $g_j$ and $r_j$ are as in Definition \ref{class}.\\

Let ${\mathcal{A}}  =(\rho ,s,s', \left\{V_j\right\}_{j=1}^s,
\left\{{r}_j\right\}_{j=1}^s)$ be an atlas in ${\mathbb{R}}^N$.
 We consider a partition of unity
$\{\psi_j\}_{j=1}^{s}$ such that $\psi_j\in C^{\infty}_c(\mathbb{R}^N)$, ${\rm supp}\,\psi_j\subset {(V_j)}_{\frac{3}{4}\rho}$, $0\le \psi_j (x) \le 1$,
 $|\nabla \psi_j(x)|\le G$
for all $x\in {\mathbb{R}}^N$
and  $j=1,\dots , s$,
where $G>0$ depends only on
${\mathcal{A}} $,
and such that $\sum_{j=1}^s\psi_j(x)=1$ for all $x\in \cup_{j=1}^s(V_j)_{\rho}$.

For $\epsilon \geq 0$ we consider  the following transformation
\begin{equation}
\label{budaout1}
T_{\epsilon}(x)=x-\epsilon\sum_{j=1}^s\xi_j\psi_j(x)\, ,\ \ \ x\in\mathbb{R}^N,
\end{equation}
where
$\xi_j={r}_j^{(-1)}((0,\dots ,1))$, which was introduced in  Burenkov and Davies~\cite{buda}.

Then we have the following variant of Lemma 18 in \cite{buda}.

\begin{lem}
\label{budaout}
Let ${\mathcal{A}}  =(\rho ,s,s', \left\{V_j\right\}_{j=1}^s,
\left\{{r}_j\right\}_{j=1}^s)$ be an atlas in ${\mathbb{R}}^N$.

Then there exist $A_1, A_2, E_1>0$ depending only on $N$ and ${\mathcal{A}}  $
 such that
\begin{equation}
\label{buda1b} \max_{0\le |\alpha |\le m}\bigl\| D^{\alpha} ( T_{\epsilon }-{\rm Id})
\bigr\|_{L^{\infty}(\mathbb{R}^N)}\le A_1\epsilon ,
\end{equation}
and such that
\begin{equation}
\label{buda1,5b} \frac{1}{2}\le 1-A_2\epsilon \le \det \nabla
T_{\epsilon }\le 1+A_2\epsilon ,
\end{equation}
for all $0\le \epsilon <E_1$. Furthermore,
\begin{equation}
\label{budaincl}
T_{\epsilon }(\Omega_1)\subset \Omega_2
\end{equation}
for all $0<\epsilon < E_1 $,
for all $\Omega_1,\Omega_2\in C({\mathcal{A}}) $ such that $\Omega_2\subset \Omega_1$ and
\begin{equation}
\label{budainclcon}
d_{{\mathcal{A}}}(\Omega_1,\Omega_2)<\frac{\epsilon }{s}.
\end{equation}
\end{lem}

{\bf Proof. } Inequalities (\ref{buda1b}), (\ref{buda1,5b}) are obvious. We now prove inclusion (\ref{budaincl}).
Let $\Omega_1, \Omega_2\in C({\mathcal{A}}) $ satisfy $\Omega_2\subset \Omega_1$ and (\ref{budainclcon}).
For all $j=1, \dots ,s$ we denote by $g_{1j}$, $g_{2j}$ respectively, the functions describing the boundaries of $\Omega_1$, $\Omega_2$  respectively, as in Definition~\ref{class} $(iii)$.  For all $x\in \cup_{j=1}^sV_j$ we set $J(x)=\{j\in \{1, \dots , s\}:\ x\in (V_j)_{\frac{3}{4}\rho} \}$.  Let $x\in \Omega_1$.  By the proof of \cite[Lemma~18]{buda} it follows that if $0<\epsilon <\frac{\rho}{4}$
then $T_{\epsilon }(x)\in \Omega_1\cap (V_j)_{\frac{\rho}{2}}$ and $d_j( T_{\epsilon }(x),  \partial \Omega_1 )\geq \epsilon \psi_j(x)$ for all $j\in J(x)$,
where $\{\psi_j\}_{j=1}^s$ is the appropriate partition of unity satisfying ${\rm supp }\,\psi_j\subset (V_j)_{\frac{3}{4}\rho}$.
Therefore
$$ \sum_{j\in J(x)} d_j( T_{\epsilon }(x), \partial \Omega_1 )\geq \epsilon \sum_{j\in J(x)} \psi_j(x)=\epsilon \sum_{j=1}^s \psi_j(x)  = \epsilon .$$ Hence there exists $\tilde \jmath\in J(x)$ such that
$
d_{\tilde \jmath}(T_{\epsilon }(x),  \partial \Omega_1  )\geq \frac{\epsilon }{s},
$
which implies
\begin{equation}
\label{budaoutpr}
(r_{\tilde \jmath}(T_{\epsilon }(x)))_N< g_{2\tilde\jmath}(\overline{r_{\tilde \jmath}(T_{\epsilon }(x)) } ).
\end{equation}
Indeed, assume to the contrary that $(r_{\tilde \jmath}(T_{\epsilon }(x)))_N\geq g_{2\tilde\jmath}(\overline{r_{\tilde \jmath}(T_{\epsilon }(x)) } )$. Then we would have
\begin{eqnarray}\lefteqn{
d_{{\mathcal{A}}}(\Omega_1,\Omega_2)\geq g_{1\tilde\jmath}(\overline { r_{\tilde\jmath} (T_{\epsilon }(x)) } )- g_{2\tilde\jmath}(\overline { r_{\tilde\jmath} (T_{\epsilon }(x)) } )}\nonumber \\ & & \qquad\quad\quad
=
g_{1\tilde\jmath}(\overline { r_{\tilde\jmath} (T_{\epsilon }(x))}  ) -(r_{\tilde \jmath} (T_{\epsilon } (x) )  )_N+(r_{\tilde \jmath} (T_{\epsilon } (x) )  )_N - g_{2\tilde\jmath}(\overline { r_{\tilde\jmath} (T_{\epsilon }(x))}  )\nonumber \\ & & \qquad\quad\quad
\geq g_{1\tilde\jmath}(\overline { r_{\tilde\jmath} (T_{\epsilon }(x) ) } ) -(r_{\tilde \jmath} (T_{\epsilon } (x) )  )_N= d_{\tilde \jmath}(T_{\epsilon }(x), \partial \Omega_1)\geq \frac{\epsilon }{s}
\end{eqnarray}
which contradicts (\ref{budainclcon}). Thus, (\ref{budaoutpr}) holds hence $ T_{\epsilon }(x)\in \Omega_2 $. \hfill $\Box $\\

\begin{thm}\label{dirthm}
Let ${\mathcal{A}}$ be an atlas in ${\mathbb{R}}^N$. Let $m\in {\mathbb{N}}$, $L,\theta >0$ and, for all $\alpha ,\beta \in {\mathbb{N}}_0^N$ with $|\alpha |=|\beta |=m$, let $A_{\alpha \beta }\in C^{0,1}(\cup_{j=1}^sV_j)$  satisfy $A_{\alpha \beta }=A_{ \beta \alpha }$, $\| A_{\alpha\beta} \|_{ C^{0,1}(\cup_{j=1}^sV_j)}$ $\le L $ and condition (\ref{elp}).

 Then for each $n\in {\mathbb{N}}$ there exist $c_n, \epsilon_n >0$ depending only on
$ n, N, {\mathcal{A}}, m, L, \theta  $   such that
\begin{equation}
\label{dirthm1}
| \lambda_{n, {\mathcal{D}}}[\Omega_1]-\lambda_{n  , {\mathcal{D}}}[\Omega_2] |\le c_n d_{\mathcal{A}}(\Omega_1,\Omega_2),
\end{equation}
 for all $\Omega_1, \Omega_2\in C({\mathcal{A}})$ satisfying $d_{\mathcal{A}}(\Omega_1,\Omega_2) <\epsilon_n $.
\end{thm}

{\bf Proof.}
Let $0<\epsilon <E_1$ where $E_1>0$ is as in Lemma~\ref{budaout}, and let
$\Omega_1,\Omega_2\in C({\mathcal{A}}) $ satisfy (\ref{budainclcon}).
We set  $\Omega_3=\Omega_1\cap\Omega_2$. Clearly, $\Omega_3\in C({\mathcal{A}}) $ and
$d_{{\mathcal{A}}}(\Omega_3,\Omega_1)$, $ d_{{\mathcal{A}}}(\Omega_3,\Omega_2)< \epsilon /s $.
By Lemma~\ref{budaout} applied to the couples of open sets $\Omega_i, \Omega_3 $ it follows that  $T_{\epsilon}(\Omega_i) \subset \Omega_3$,
$i=1,2$.
By the monotonicity of the eigenvalues with respect to inclusion it follows that
\begin{equation}
\label{dirproof1}
\lambda_{n, {\mathcal{D}}}[\Omega_i]\le \lambda_{n, {\mathcal{D}}}[\Omega_3]\le \lambda_{n, {\mathcal{D}}}[T_{\epsilon }(\Omega_i)],\ \ \ i=1,2.
\end{equation}
Since in Lemma~\ref{unifbound} $\Lambda_n$ depends only on $n,N,{\mathcal{A}},m$ and $L$, in Corollary~\ref{diffeocorol} $c$ depends only
on $ N, {\mathcal{A}}, m, B_1, B_2, L  $ and $\theta$, and in Lemma~\ref{budaout} $E_1$ and $A_1$ depend only on $N$ and ${\mathcal{A}}$,
by (\ref{diffeocorol1}), (\ref{newnumb}) and (\ref{buda1b}) it follows that there exist $\tilde c_n, \tilde \epsilon_n>0$ such that
$$
\lambda_{n,{\mathcal{D}}}[\Omega_3]-\lambda_{n,{\mathcal{D}}}[\Omega_i]\le \lambda_{n,{\mathcal{D}}}[T_{\epsilon }(\Omega_i)]
-\lambda_{n,{\mathcal{D}}}[\Omega_i]\le \tilde c_n \epsilon,\ \ \ i=1,2,
$$
if $0<\epsilon <\tilde \epsilon_n$.
Hence
$$
|\lambda_{n,{\mathcal{D}}}[\Omega_1]-\lambda_{n,{\mathcal{D}}}[\Omega_2]|\le\max_{i=1,2}\{    \lambda_{n,{\mathcal{D}}}[\Omega_3]-
\lambda_{n,{\mathcal{D}}}[\Omega_i]    \}\le \tilde c_n\epsilon .
$$
Take here $ \epsilon=2s d_{{\mathcal{A}}}(\Omega_1,\Omega_2) $, then inequality (\ref{dirthm1}) holds with $c_n=2s\tilde c_n$
if $d_{{\mathcal{A}}}(\Omega_1,\Omega_2)<\epsilon_n=\tilde\epsilon_n/(2s)$.\hfill $\Box$\\

Let
${\mathcal{A}}= (  \rho , s,s', \{V_j\}_{j=1}^s, \{r_j\}_{j=1}^{s} ) $ be an atlas in ${\mathbb{R}}^N$.
For all $x\in V'=\cup_{j=1}^{s'}V_j$ we set $J'(x)=\{j\in \{1, \dots , s'\}:\ x\in V_j  \}$.
Let $\Omega\in C({\mathcal{A}})$.
Then
we set
$$
d_{\mathcal{A}}(x, \partial \Omega )= \max_{j\in J'(x)}d_j(x, \partial \Omega ),
$$
for all $x\in V'$,
where $d_j(x, \partial \Omega )$ is defined in (\ref{dj}). Observe that if $\Omega\in C({\mathcal{A}})$ then $\partial \Omega \subset V'$.
Therefore if $\Omega_1, \Omega_2\in C({\mathcal{A}})$ then
\begin{equation}
\label{atl}
d_{\mathcal{A}}( \Omega_1,  \Omega_2 )=\sup_{x\in \partial \Omega_1}d_{\mathcal{A}}(x, \partial \Omega_2 )=
\sup_{x\in \partial \Omega_2}d_{\mathcal{A}}(x, \partial \Omega_1 ).
\end{equation}

For all $\epsilon >0$ we set
$$
\Omega_{\epsilon , {\mathcal{A}}}=  \Omega \setminus \{x\in  V':\ d_{\mathcal{A}}(x, \partial \Omega )\le \epsilon  \},
$$
$$
\Omega^{\epsilon , {\mathcal{A}}}=  \Omega \cup \{x\in  V':\ d_{\mathcal{A}}(x, \partial \Omega )< \epsilon  \}.
$$

\begin{lem}\label{cinquesedici}
Let
${\mathcal{A}}$ be an atlas in ${\mathbb{R}}^N$ and $\epsilon >0$.
If $\Omega_1$ and $\Omega_2$ are
two open sets in $C({\mathcal{A}})$ satisfying the inclusions
\begin{equation} \label{maininclA}
(\Omega_1)_{\epsilon , {\mathcal{A}} }\subset \Omega_2 \subset (\Omega_1)^{\epsilon , {\mathcal{A}} }
\end{equation}
or
\begin{equation} \label{maininclbisA}
(\Omega_2)_{\epsilon , {\mathcal{A}}}\subset \Omega_1 \subset (\Omega_2)^{\epsilon , {\mathcal{A}}},
\end{equation}
 then
\begin{equation}
\label{emibisA}
d_{{\mathcal{A}}}( \Omega_1, \Omega_2)\le \epsilon .
\end{equation}
\end{lem}

{\bf Proof.}  Assume that inclusion (\ref{maininclA}) holds. Let $x\in \partial \Omega_2$. We consider three cases. Case $x\in \Omega_1$. Since $x\notin \Omega_2$ then $x\notin (\Omega_1)_{\epsilon , {\mathcal{A}}}$ hence by definition of $(\Omega_1)_{\epsilon  , {\mathcal{A}}}$ it follows that $x\in V'$ and
$d_{\mathcal{A}}(x, \partial \Omega_1)\le \epsilon $. Case $x\in \partial \Omega_1$. Obviously $d_{\mathcal{A}}(x, \partial \Omega_1)=0$. Case $x\notin \overline{ \Omega}_1$. In this case there exists
a sequence $x_n\in \Omega_2\setminus \overline{\Omega}_1, \ n\in {\mathbb{N}}$ converging to $x$. Since $x_n\notin \overline{\Omega}_1$ then $d_{{\mathcal{A}}}(x_n ,\partial \Omega_1)<\epsilon $ because $x_n \in (\Omega_1)^{\epsilon , {\mathcal{A}}}$. By observing that $J'(x_n)=J'(x)$ for all
$n$ sufficiently large,
one can pass to the limit and obtain
$d_{{\mathcal{A}}}(x ,\partial \Omega_1)\le \epsilon $. Thus, in any case we have that  $d_{{\mathcal{A}}}(x ,\partial \Omega_1)\le \epsilon $
for all $x\in \partial \Omega_2 $ and (\ref{emibisA}) follows by (\ref{atl}). The same argument applies when  inclusion (\ref{maininclbisA}) holds.
\hfill $\Box$\\

\begin{corol}
\label{lastcorADIR}
Let ${\mathcal{A}}$ be an atlas in ${\mathbb{R}}^N$. Let $m\in {\mathbb{N}}$,  $L, \theta >0$ and, for all $\alpha ,\beta \in {\mathbb{N}}_0^N$ with $|\alpha |=|\beta |=m$, let $A_{\alpha \beta }\in C^{0,1}(\cup_{j=1}^sV_j)$  satisfy $A_{\alpha \beta }=A_{ \beta \alpha }$, $\| A_{\alpha\beta} \|_{ C^{0,1}(\cup_{j=1}^sV_j)}\le L $ and condition (\ref{elp}).

 Then for each $n\in {\mathbb{N}}$ there exist $c_n, \epsilon_n >0$ depending only on
$n, N, {\mathcal{A}},  m,   L, \theta  $   such that
\begin{equation}
\label{maincordavstyA}
| \lambda_{n, {\mathcal{D}}}[\Omega_1]-\lambda_{n, {\mathcal{D}} }[\Omega_2] |\le c_n  \epsilon ,
\end{equation}
for all $0<\epsilon <\epsilon_n$ and
 for all $\Omega_1, \Omega_2\in C({\mathcal{A}})$ satisfying (\ref{maininclA}) or (\ref{maininclbisA}).
\end{corol}

{\bf Proof. } Inequality (\ref{maincordavstyA}) follows by inequality (\ref{dirthm1}) and  inequality (\ref{emibisA}). \hfill $\Box$\\

\section{Estimates for Neumann eigenvalues via the atlas distance}
\label{section6}

In this section we prove Theorem \ref{neuthm}. The proof is based on Lemmas~\ref{patchesthm} and \ref{newpatches}.

\begin{defn}
Let $U$ be an open set in ${\mathbb{R}}^N$ and $\rho$ a rotation. We say that $U$ is a `$\rho$-patch'
if there exist an open set $G_U\subset {\mathbb{R}}^{N-1}$ and functions $\varphi_{U}, \psi_{U}:G_U\to {\mathbb{R}}$ such that
$$
\rho (U)=\left\{(\bar x, x_N)\in {\mathbb{R}}^N:\ \psi_U(\bar x)<x_N<\varphi _U(\bar x),\ \bar x \in G_U  \right\}.
$$

The `thickness' of the $\rho$-patch is defined by
$$
R _U=\inf_{\bar x\in G_U}(\varphi_U(\bar x)-\psi_U (\bar x) );
$$
the `thinness' of the $\rho$-patch is defined by
$$
S _U=\sup_{\bar x\in G_U}(\varphi_U(\bar x)-\psi_U (\bar x) ).
$$

\end{defn}

If $\Omega_2\subset \Omega_1$ and $\Omega_1\setminus \Omega_2$ is covered by a finite number of $\rho$-patches contained in $\Omega_1$, then
we can estimate $\lambda_{n,{\mathcal{N}}}[\Omega_2]- \lambda_{n,{\mathcal{N}}}[\Omega_1]$ via the thinness of the patches.

\begin{lem}
\label{patchesthm}
Let $m\in {\mathbb{N}}$ and $\Omega_1$ be an open set in ${\mathbb{R}}^N$ such that the embedding $W^{m,2}(\Omega_1 )\subset W^{m-1,2}(\Omega_1 )$ is compact.
For all $\alpha ,\beta \in
 \mathbb{N}_0^{N} $
with $|\alpha |=|\beta |=m$, let
$A_{\alpha \beta }$ be bounded measurable real-valued
functions defined on $\Omega_1$,  satisfying
$A_{\alpha \beta }= A_{ \beta\alpha  }$ and condition (\ref{elp}) in $\Omega_1$.
Let $\sigma \in {\mathbb{N}}$, $R>0$.

Assume that  $\Omega_2\subset \Omega_1$ is such that the embedding $W^{m,2}(\Omega_2 )\subset W^{m-1,2}(\Omega_2 )$ is compact
and there exist rotations $\{\rho_j\}_{j=1}^{\sigma }$ and two sets
 $\{ U_j\}_{j=1}^{\sigma }$, $\{ \tilde U_j\}_{j=1}^{\sigma }$ of $\rho_j$-patches $U_j$ and $\tilde U_j$ satisfying the following properties
\begin{itemize}
\item[(a)] $U_j\subset \tilde U_j \subset \Omega_1$,
for all $j=1, \dots , \sigma $;
\item[(b)]  $G_{U_j}=G_{\tilde U_j}$, $\varphi_{U_j}=\varphi_{\tilde U_j}$, for all $j=1, \dots , \sigma $;
\item[(c)]  $R_{\tilde U_j}> R$, for all $j=1, \dots , \sigma $;
\item[(d)]  $\Omega_1\setminus \Omega_2\subset \cup_{j=1}^{\sigma }U_j$.
\end{itemize}

Then there exists $d>0$ depending only on $N, m,R$ such that for all $n\in {\mathbb{N}}$
\begin{equation}
\label{patchest}
\lambda_{n,{\mathcal{N}}}[\Omega _2]\le \lambda_{n,{\mathcal{N}}}[\Omega _1](1+d_n\max _{j=1, \dots , \sigma }S_{U_j}),
\end{equation}
if
$\max_{j=1, \dots , \sigma }S_{U_j}<d_n^{-1}$,
where
\begin{equation}
\label{patchesconst}d_n= 2\sigma d (1+\theta^{-1}\lambda_{n, {\mathcal{N}}}[\Omega _1]) .
\end{equation}
\end{lem}

{\bf Proof.} By $(a)$ and $(b)$ it follows  that  $\psi_{\tilde U_j}\le \psi_{ U_j}$ for all $j=1, \dots , \sigma $.
Let $u\in W^{m,2}(\Omega_1)$. By $(d)$
\begin{equation}
\label{p1}
\int_{\Omega_1\setminus \Omega_2}|u|^2dy\le \sum_{j=1}^{\sigma }\int_{U_j}|u|^2dy= \sum_{j=1}^{\sigma }\int_{\rho_{j}(U_j)}|u\circ \rho_{j}^{(-1)}|^2dx.
\end{equation}
Let, for brevity,  $v_j=u\circ \rho_{j}^{(-1)}$.  Clearly,
\begin{equation}
\label{p2}
\int_{\rho_{j}(U_j)}|u\circ \rho_{j}^{-1}|^2=\int_{G_{U_j}}\int_{\psi_{U_j}(\bar x)}^{\varphi_{U_j}(\bar x)}|v_j(\bar x, x_N)|^2d\bar xdx_N.
\end{equation}

Since $v_j\in W^{m,2}(\rho_{j}(\tilde U_j))$ it follows that for almost all $\bar x\in G_{\tilde U_j}$ the function $v_j(\bar x, \cdot )$ belongs to the
space $W^{m,2}( \psi_{\tilde U_j}(\bar x),  \varphi_{ \tilde U_j}(\bar x)  ) $. Moreover, by (c) it follows that $ \varphi_{\tilde U_j}(\bar x)-\psi_{\tilde U_j}(\bar x)\geq R   $. Thus by Burenkov~\cite[Thm.~2,~p.127]{bu} there exists $\tilde d>0$ depending only on $m, R$ such that
\begin{eqnarray}\label{bur}\lefteqn{
\| v_j(\bar x, \cdot )  \|^2_{L^{\infty }(\psi_{\tilde U_j}(\bar x), \varphi_{\tilde U_j}(\bar x) ) }\le \tilde d  \biggl(\|  v_j(\bar x ,\cdot  )\|^2_{L^{2}( \psi_{\tilde U_j}(\bar x), \varphi_{\tilde U_j}(\bar x) )}   \biggr. }\nonumber  \\ & & \qquad\qquad\qquad\qquad\qquad\qquad\qquad +  \biggl.   \biggl\| \frac{\partial^m v_j}{\partial x_N^m}(\bar x ,\cdot  )\biggr\|^2_{L^{2}( \psi_{\tilde U_j}(\bar x), \varphi_{\tilde U_j}(\bar x) )}\biggr).
\end{eqnarray}
By inequality (\ref{bur}) and property $(b)$
\begin{eqnarray}
\label{p3}
\lefteqn{
\int_{G_{U_j}}\int_{\psi_{U_j}(\bar x)}^{\varphi_{U_j}(\bar x)}|v_j(\bar x, x_N)|^2d\bar xdx_N} \nonumber \\
& & \qquad \le  \int_{G_{U_j}} ( \varphi_{U_j}(\bar x)-\psi_{U_j}(\bar x)  )\| v_j(\bar x , \cdot ) \|^2_{L^{\infty }(  \psi_{\tilde U_i}(\bar x),  \varphi_{ U_i}(\bar x) ) }d\bar x \nonumber \\ & & \qquad
\le \tilde d S_{U_j} \biggl( \| v_j \|^2_{L^{2}(\rho_{j}(\tilde U_j)  )}+ \left\| \frac{\partial^mv_j}{\partial x_N^m}    \right\|^2_{L^{2}(\rho_{j}(\tilde U_j)   )} \biggr)\nonumber \\  & & \qquad \le
 d  S_{U_j}\biggl(  \| u \|^2_{L^{2}(\Omega _1)}+\sum_{|\alpha |=m} \| D^{\alpha }u \|^2_{L^{2}(\Omega_1)}  \biggr),
\end{eqnarray}
where  $d>0$ depends only on $N, m, R$.

Let $\psi_n[\Omega_1] $, $n\in {\mathbb{N}}$, be an orthonormal sequence of eigenfunctions corresponding to the
eigenvalues $\lambda_{n,{\mathcal{N}}}[\Omega_1]$. We denote by $L_n[\Omega_1]$ the linear subspace of $W^{m,2}(\Omega_1)$ generated by the
$\psi_1[\Omega_1], \dots , \psi_n[\Omega _1]$.
If $u\in L_{n}[\Omega_1]$ and $\| u\|_{L^2(\Omega_1)}=1$ then by (\ref{p1}), (\ref{p2}
), (\ref{p3})
\begin{eqnarray}
\lefteqn{ \int_{\Omega_1\setminus\Omega_2} |u|^2 \le
\sigma d \max_{j=1,\dots , \sigma } S_{U_j}(1+\theta^{-1}Q_{\Omega_1}(u)  )}\nonumber \\ & & \qquad\qquad\qquad\qquad  \le \sigma d \max_{j=1,\dots , \sigma } S_{U_j}
(1+\theta^{-1}\lambda_n[\Omega_1] ).
\end{eqnarray}

Let ${\mathcal{T}}_{12}$ be the restriction operator from $\Omega_1$ to $\Omega_2$. Clearly, ${\mathcal{T}}  _{12}$ maps $W^{m,2}(\Omega_1)$
to $W^{m,2}(\Omega_2)$. For all $n\in {\mathbb{N}}$ and for all $u\in L_{n}[\Omega_1]$, $\| u\|_{L^2(\Omega_1)}=1$ we have
\begin{eqnarray}
\| {\mathcal{T}}  _{12}u\|^2_{L^2(\Omega_2)}= \int_{\Omega_1}|u|^2-\int_{\Omega_1\setminus\Omega_2}|u|^2\geq 1-     \sigma d \max_{j=1,\dots , \sigma } S_{U_j}
(1+\theta^{-1}\lambda_n[\Omega_1] )
\end{eqnarray}
and
$$
Q_{\Omega_2}({\mathcal{T}}  _{12}u)\le Q_{\Omega_1 }(u)\le \lambda_n[\Omega_1]
$$
because $ \sum_{|\alpha |=|\beta |=m}A_{\alpha \beta}\xi_{\alpha}\bar\xi_{\beta} \geq 0$ for all $\xi_{\alpha },\xi_{\beta}\in {\mathbb{C}}$.
Thus, in the terminology of \cite{bula},  $ {\mathcal{T}} _{12}$ is a transition operator from $H_{W^{m,2}(\Omega_1)}$ to $H_{W^{m,2}(\Omega_2)}$ with the measure of vicinity $\delta (H_{W^{m,2}(\Omega_1)},$ $ H_{W^{m,2}(\Omega_2)})=\max_{j=1,\dots , \sigma } S_{U_j}$ and the parameters $a_n=  \sigma d  (1+\theta^{-1}\lambda_n[\Omega_1] ) $, $b_n=0$. Thus,
by the general spectral stability theorem \cite[Theorem~3.2]{bula} it follows that
$$
\lambda_{n,{\mathcal {N}}}[\Omega_2 ]\le \lambda_{n,{\mathcal {N}}}[\Omega_1 ]+2(a_n\lambda_{n,{\mathcal{N}}}[\Omega _ 1]+b_n)\delta
(H_{W^{m,2}(\Omega_1)}, H_{W^{m,2}(\Omega_2)})$$
if $\delta (H_{W^{m,2}(\Omega_1)}, H_{W^{m,2}(\Omega_2)})< (2a_n)^{-1}$ which immediately gives (\ref{patchest}).\hfill $\Box$\\

\begin{lem}
\label{semi}
Let ${\mathcal{A}}$ be an atlas in ${\mathbb{R}}^N$. Let $m\in {\mathbb{N}}$,  $L,\theta >0$ and, for all $\alpha ,\beta \in {\mathbb{N}}_0^N$ with $|\alpha |=|\beta |=m$, let $A_{\alpha \beta }\in L^{\infty }(\cup_{j=1}^sV_j)$  satisfy $A_{\alpha \beta }=A_{ \beta \alpha }$, $\| A_{\alpha\beta} \|_{ L^{\infty }(\cup_{j=1}^sV_j)}\le L $ and condition (\ref{elp}).
Then for each $n\in {\mathbb{N}}$ there exist $c_n,\epsilon_n>0 $ depending only on $n, N, {\mathcal{A}} , m, L, \theta $, such that
\begin{equation}
\label{semiest}
\lambda_{n,{\mathcal{N}}}[\Omega_2]\le \lambda_{n,{\mathcal{N}}}[\Omega_1]+c_n d_{ {\mathcal{A}} }(\Omega_1,\Omega_2),
\end{equation}
for all $\Omega_1, \Omega _2\in C({\mathcal{A}})$ satisfying $\Omega_2\subset \Omega_1$ and
$d_{ {\mathcal{A}} }(\Omega_1,\Omega_2 )<\epsilon_n$.
\end{lem}

{\bf Proof.} Let  ${\mathcal{A}}= (  \rho , s,s', \{V_j\}_{j=1}^s, \{r_j\}_{j=1}^{s} ) $,
 $\Omega_1, \Omega_2\in C({\mathcal{A}})$ and $\Omega_2\subset \Omega_1$.
For all $j=1, \dots ,s$ we denote by $g_{1j}$, $g_{2j}$ respectively, the functions describing the boundaries of $\Omega_1$, $\Omega_2$  respectively, as in Definition~\ref{class}  $(iii)$.
We consider two sets of  $\{ U_j \}_{j=1}^s$
$\{\tilde  U_j \}_{j=1}^s$  of $r_j$-patches $U_j$, $\tilde U_j$ defined as follows:
$$
\tilde U_j = r_j^{(-1)}(\{ (\bar x , x_N):\ \bar x\in W_j,\ a_{Nj}<x_N <g_{1j}(\bar x )  \}   ) ,
$$
$$
  U_j = r_j^{(-1)}(\{ (\bar x , x_N):\ \bar x\in W_j,\   g_{2j}(\bar x )<x_N <g_{1j}(\bar x )  \},
$$
where $W_j$ and $a_{Nj}$ are as in Definition~\ref{class}.
Observe that conditions $(a), (b), (c), (d)$ of Lemma \ref{patchesthm} are satisfied with $\sigma =s$ and $R=\rho$. Moreover, $\max_{j=1, \dots , \sigma }S_{U_j}=d_{{\mathcal{A}}}(\Omega_1,\Omega_2)$. Thus by applying Lemma~\ref{patchesthm}  to the open sets $\Omega_1$, $\Omega_2$ and the sets of patches defined above,
and by Lemma~\ref{unifbound} we immediately deduce the validity of (\ref{semiest}).\hfill $\Box $\\

Our next aim is to consider the case $\Omega_2=T_{\epsilon }(\Omega _1)$ where
 $T_{\epsilon}$ is the map defined in (\ref{budaout1}).

\begin{lem}
\label{newpatches}
Let ${\mathcal{A}}$ be an atlas in ${\mathbb{R}}^N$.
Then there exist $\epsilon_0, A, R >0$ and $\sigma\in \mathbb{N}$ depending only on $N, {\mathcal{A}}  $, and for each open set $\Omega\in C({\mathcal{A}}  )$ and for each
$0<\epsilon <\epsilon_0$
there exist
rotations $\{\rho_j\}_{j=1}^{\sigma }$ and
sets $\{U_j\}_{j=1}^{\sigma} $, $ \{\tilde U_j\}_{j=1}^{\sigma} $ of $\rho_j$-patches $U_j$, $\tilde U_j$
satisfying conditions $(a)$, $(b)$, $(c)$, $(d)$ in Lemma \ref{patchesthm} with  $\Omega_1=\Omega$ and $\Omega _2=T_{\epsilon }(\Omega )$ and such that
$\max_{j=1, \dots , \sigma}S_{U_j}<A\epsilon  $.
\end{lem}

{\bf Proof.} In fact, we shall prove that there exist a family  of rotations $ \{\rho_j\}_{j=1}^{\sigma }   $,
a family $ \{G_j\}_{j=1}^{\sigma } $ of bounded open sets in $\mathbb{R}^{N-1}$, and a family  $ \{\varphi_j\}_{j=1}^{\sigma } $ of functions $\varphi_j$ continuous on $\overline{ G}_j$  such that for all $0<\epsilon < \epsilon_0$
\begin{equation}
\label{strip}
\Omega\setminus T_{\epsilon}(\Omega)\subset \bigcup\limits_{j=1}^{\sigma} {U}_j^{[\epsilon]}
\end{equation}
and
\begin{equation}
\label{strip1}
{U}_j^{[\epsilon]}\subset \tilde {U}_j\subset \Omega ,
\end{equation}
where the $\rho_j$-patches ${U}_j^{[\epsilon]}  $, $\tilde {U}_j$ are defined by
\begin{equation}
\label{strip2}
\rho_j( {U}_j^{[\epsilon]})=\left\{(\bar x , x_N)\in \mathbb{R}^{N-1}:\ \varphi_j(\bar x)-A\epsilon< x_N <\varphi_j(\bar x),\ \bar x\in G_j \right\},
\end{equation}
and
\begin{equation}
\label{strip3}
\rho_j(\tilde {U}_j)=\left\{(\bar x , x_N)\in \mathbb{R}^{N-1}:\ \varphi_j(\bar x)-R< x_N <\varphi_j(\bar x),\ \bar x\in G_j \right\}.
\end{equation}

Let ${\mathcal{A}}  =(\rho ,s,s', \left\{V_j\right\}_{j=1}^s,
\left\{{r}_j\right\}_{j=1}^s)$ and $\Omega \in C({\mathcal{A}})$. We split the proof into four steps.

{\it Step 1.} Let for each non-empty set $J\subset \{1, \dots , s'\}$, $V_J=\cap_{j\in J}(V_j)_{\frac{\rho}{2}}$ and $d={\rm dim\, Span}\{\xi_{j} \}_{j\in J}$. Recall that $\xi_j=r_j^{(-1)}(0, \dots , 1)$. By the proof of Lemma 19 in \cite{buda} it follows that there exist vectors $\xi_J\equiv \xi_{J1}, \xi_{J2}, \dots ,\xi_{Jd}$ and a
rotation $r_J$ such that:

1) $\xi_J, \xi_{J2}, \dots ,\xi_{Jd}$ is an orthonormal basis for  ${\rm Span}\{\xi_{j} \}_{j\in J}$ and $r_J(\xi_J)=e_N$, $r_J(\xi_{J2})=e_{N-1}, \dots ,r_J(\xi_{Jd})=e_{N-d+1}$,

2) there exist continuous functions $\varphi_J, \psi_J$ defined on $\overline{ G}_J$ where $G_J={\rm Pr}_{x_N=0}$ $r_J(V_J\cap \Omega )$  (${\rm Pr}_{x_N=0}$ denotes the orthogonal projector onto the hyperplane with the equation $x_N=0$) such that
\begin{equation}
\label{parstrip}
r_J(V_J\cap\Omega )=\left\{(\bar x, x_N)\in {\mathbb {R}}^{N-1}:\ \psi_J(\bar x)<x_N<\varphi_J(\bar x) \right\}
\end{equation}
and such that
\begin{equation}
\label{parstrip1}
\left\{(\bar x, y)\in {\mathbb {R}}^{N-1}:\ y<\varphi_J(\bar x),\ \bar x\in G_J,\ (\bar x, y)\in r_J(V_j),\ \forall\, j\in J \right\}\subset r_J(\Omega ),
\end{equation}

3) the function $\varphi_J$ satisfies the Lipschitz condition with respect to the variables $v=(x_{N-d+1}, \dots , x_{N-1})$ uniformly with respect
to the va\-ria\-bles $u=(x_1, \dots , x_{N-d})$ on $G_J$, {\it i.e.,}
$$
|\varphi_J(u,v)-\varphi_J(u,w)|\le L_J|v-w|,\ \ \ \ {\rm for\ all}\ (u,v),(u,w)\in G_J,
$$
where $L_J>0$ depends only on $\{V_j \}_{j=1}^s$ and $\{r_j \}_{j=1}^s$.

Observe that by (\ref{parstrip}) and (\ref{parstrip1}) it follows that
\begin{equation}
\label{parstrip2}
\left\{(\bar x, x_N)\in {\mathbb {R}}^{N-1}:\ \bar x\in G_J,\ \psi_J(\bar x)-\frac{\rho }{4}<x_N<\varphi_J(\bar x) \right\}\subset r_J(\Omega )
\end{equation}
because the distance of $(\bar x, \psi_J(\bar x))$ to the boundary of $r_J(V_j)$ is greater than $\frac{\rho}{2}$ for all $j\in J$ and $\bar x\in G_J$.

{\it Step 2.} For $x\in \mathbb{R}^{N}$ let as in \cite{buda}
$$
J(x)= \{j\in\{1, \dots ,s \}:\ x\in (V_j)_{\frac{3 }{4 }\rho} \}.
$$
Observe that $J(x)\subset \{1, \dots , s'\}$ if $x\in \partial \Omega$.
The inclusion ${\rm supp}\psi_j\subset (V_j)_{\frac{3 }{4 }\rho }$ implies that $\psi_j(x)=0$ for $j\notin J(x)$ and
$$
T_{\epsilon }(x)=x-\epsilon \sum_{j\in J(x)}\xi_j\psi_j(x),\ \ \ x\in \mathbb{R}^N.
$$

For any  subset $J\subset \{1, \dots , s \}$ we set
$$
\tilde V_J=\left\{ x\in \mathbb{R}^N:\ J(x)=J\right\}
$$
so that $\displaystyle {\mathbb{R}}^N= \mathring{\cup}_{\substack{J\subset \{1, \dots , s'\}}}\tilde V_J$ and
$$
T_{\epsilon}(x)=x-\epsilon \sum_{j\in J}\xi_j\psi_j(x),\ \ \ {x\in \tilde V_J}.
$$

{\it Step 3.} Let $x\in \tilde V_J\cap \partial \Omega$. Since $\| T_{\epsilon }- {\rm Id}\|_{\infty }\le \epsilon$ and $T_{\epsilon }(\Omega )\subset \Omega$ we have that $T_{\epsilon }(x) \in V_J\cap \Omega $ for all $ 0<\epsilon \le \frac{\rho}{4}$.  Let
$r_J(x)=(\beta^{(1)}, \beta^{(2)}, \beta_N)$ where $\beta^{(1)}=(\beta_1, \dots , \beta_{N-d})$, $\beta^{(2)}=(\beta_{N-d+1}, \dots , \beta_{N-1})$ and $\beta_N=\varphi_J(\beta^{(1)}, \beta^{(2)}) $. Since $T_{\epsilon }(x)-x\in {\rm Span} \{\xi_j\}_{j\in J}$, it follows that $r_J(T_{\epsilon }(x))= (\beta^{(1)}, \gamma^{(2)}, \gamma_N)$ for some $\gamma^{(2)}=(\gamma_{N-d+1}, \dots , \gamma_{N-1})$ and $\gamma_N$.
Since $T_{\epsilon }(x) \in V_J\cap \Omega $, for the distance $d_J(T_{\epsilon }(x) ) $ of $T_{\epsilon }(x)$ from $\partial \Omega $ in the direction of the vector $\xi_J$ we have
\begin{eqnarray}
\lefteqn{ d_J(T_{\epsilon }(x) )=\varphi_J(\beta ^{(1) }, \gamma ^{(2)})-\gamma_N}\nonumber\\ & & \qquad\qquad\qquad
=\varphi_J(\beta ^{(1) }, \gamma ^{(2)})-\beta_N+\beta_N-\gamma_N\nonumber\\
& & \qquad\qquad\qquad=\varphi_J(\beta ^{(1) }, \gamma ^{(2)})-\varphi_J(\beta ^{(1) }, \beta ^{(2)})+\beta_N-\gamma_N\nonumber\\
& & \qquad\qquad\qquad\le L_J |\gamma^{(2)}-\beta^{(2)} |+ |\gamma_N-\beta_N|\nonumber\\
& & \qquad\qquad\qquad\le (L_J+1)|r_J(T_{\epsilon}(x))-r_J(x)|\nonumber\\
& & \qquad\qquad\qquad= (L_J+1)|T_{\epsilon}(x)-x| .\label{dJ}
\end{eqnarray}
Let
\begin{equation}
\label{A}
A=\max_{\substack{J\subset \{1, \dots , s' \} \\ J\ne\emptyset  }}(L_J+1)
\end{equation}
and $ {U}_J^{[\epsilon]} $ be defined by
\begin{equation}
\label{strip2bis}
\rho_J( {U}_J^{[\epsilon]})=\left\{(\bar x , x_N)\in \mathbb{R}^{N-1}:\ \varphi_J(\bar x)-A\epsilon< x_N <\varphi_J(\bar x),\ \bar x\in G_J \right\}.
\end{equation}
Then by (\ref{dJ}) it follows that $T_{\epsilon }(\tilde V_J\cap\partial \Omega  )\subset {U}_J^{[\epsilon]}$ and
\begin{equation}
\label{covering}
T_{\epsilon}(\partial \Omega )= \bigcup_{\substack{J\subset \{1, \dots , s' \} \\ J\ne\emptyset  }}T_{\epsilon}( \tilde V_J\cap\partial \Omega  )\subset
\bigcup_{\substack{J\subset \{1, \dots , s' \} \\ J\ne\emptyset  }} {U}_J^{[\epsilon]}.
\end{equation}

{\it Step 4.} Let $y\in \Omega \setminus \bigcup_{  \substack{J\subset \{1, \dots , s' \} \\ J\ne\emptyset  }      } {U}_J^{[\epsilon]}$.
By the definition of ${U}_J^{[\epsilon]}$ it follows that $y\notin {U}_J^{[\epsilon ']}$ for all $0<\epsilon '\le \epsilon$.
Thus by (\ref{covering}) it follows that $y\notin T_{\epsilon ' }(\partial \Omega )$. This implies that the topological degree ${\rm deg }(\Omega , T_{\epsilon ' }, y)$ of
the triple $(\Omega , T_{\epsilon ' }, y)$ is well defined (see {\it e.g.}, Deimling~\cite[\S 1]{dei}) and by homotopy invariance ${\rm deg }(\Omega , T_{\epsilon ' }, y)={\rm deg }(\Omega , T_{0 }, y)=1$ for all $0<\epsilon '\le \epsilon$. Thus the equation $T_{\epsilon }(x)=y$ has a solution $x\in \Omega$ hence $y\in T_{\epsilon }(\Omega )$ (see {\it e.g.}, Deimling~\cite[Thm.~3.1]{dei}).
This shows that $\Omega \setminus T_{\epsilon }(\Omega ) \subset \bigcup_{\substack{J\subset \{1, \dots , s' \} \\ J\ne\emptyset  }} {U}_J^{[\epsilon]}$.

To complete the proof of the lemma it suffices to choose
 $\sigma$ to  be the number of nonempty subsets of $\{1, \dots , s'\}$, $\epsilon_0=\frac{\rho }{4}$, $R=\frac{\rho }{4}$ (see (\ref{parstrip2})), and $A$  as in (\ref{A}). \hfill $\Box$

\begin{thm}\label{neuthm}
Let ${\mathcal{A}}$ be an atlas in ${\mathbb{R}}^N$. Let $m\in {\mathbb{N}}$,  $L,\theta >0$ and, for all $\alpha ,\beta \in {\mathbb{N}}_0^N$ with $|\alpha |=|\beta |=m$, let $A_{\alpha \beta }\in C^{0,1}(\cup_{j=1}^sV_j)$  satisfy $A_{\alpha \beta }=A_{ \beta \alpha }$, $\| A_{\alpha\beta} \|_{ C^{0,1}(\cup_{j=1}^sV_j)}\le L $ and condition (\ref{elp}).

 Then for each $n\in {\mathbb{N}}$ there exist $c_n, \epsilon_n >0$ depending only on
$n, N,   {\mathcal{A}}, m, L, \theta  $   such that
\begin{equation}
\label{neuthm1}
| \lambda_{n, {\mathcal{N}}}[\Omega_1]-\lambda_{n  , {\mathcal{N}}}[\Omega_2] |\le c_n d_{\mathcal{A}}(\Omega_1,\Omega_2),
\end{equation}
 for all $\Omega_1, \Omega_2\in C({\mathcal{A}})$ satisfying $d_{\mathcal{A}}(\Omega_1,\Omega_2) <\epsilon_n $.
\end{thm}

{\bf Proof.}
 In this proof $c_n$, $\epsilon_n$  denote positive constants depending only on some of the parameters $n, N,  {\mathcal{A}}, m, L, \theta  $ and their value is not necessarily the same for all the inequalities below.

Let $E_1>0$ be as in Lemma~\ref{budaout}. Let $0<\epsilon <E_1$ and $\Omega_1,\Omega_2\in C({\mathcal{A}}) $  satisfy (\ref{budainclcon}). We set  $\Omega_3=\Omega_1\cap\Omega_2$. Clearly, $\Omega_3\in C({\mathcal{A}}) $ and
$d_{{\mathcal{A}}}(\Omega_3,\Omega_1), d_{{\mathcal{A}}}(\Omega_3,\Omega_2)< \frac{\epsilon}{s} $. By Lemma~\ref{budaout} applied to the couple of open sets $\Omega_1, \Omega_3 $ it follows that  $T_{\epsilon}(\Omega_1) \subset \Omega_3$ hence
\begin{equation}
\label{last1}
T_{\epsilon}(\Omega_3)\subset T_{\epsilon}(\Omega_1) \subset \Omega_3.
\end{equation}

We now apply Lemma \ref{newpatches} to the set $\Omega =\Omega_3$. It follows that if $0<\epsilon <\epsilon_0 $ there exist rotations $\{\rho_j\}_{j=1}^{\sigma }$ and two sets
 $\{U_j\}_{j=1}^{\sigma} $, $ \{\tilde U_j\}_{j=1}^{\sigma} $ of $\rho_j$-patches $U_{j}$, $\tilde U_j$
satisfying conditions $(a)$, $(b)$, $(c)$, $(d)$ in Lemma \ref{patchesthm}  with $\Omega_1$ replaced by $\Omega_3$ and $\Omega_2$ replaced by
$T_{\epsilon}(\Omega_3)$, and such that
$\max_{j=1, \dots , \sigma}S_{U_j}<A\epsilon  $. In particular,
\begin{equation}
\label{last2}
\Omega_3\setminus T_{\epsilon}(\Omega_3)\subset \cup_{j=1}^{\sigma }U_j,
\end{equation}
hence by (\ref{last1}), (\ref{last2}) it follows
\begin{equation}
\label{last3}
\Omega_3\setminus T_{\epsilon}(\Omega_1)\subset \cup_{j=1}^{\sigma }U_j.
\end{equation}
Now we apply Lemma \ref{patchesthm} to the couple  of open sets $\Omega_3$, $T_{\epsilon}(\Omega_1)$ by using the sets of patches defined above.
Since $\max_{j=1, \dots , \sigma}S_{U_j}<A\epsilon$,  by  Lemma \ref{patchesthm} it follows that
if $A\epsilon <d_n^{-1} $ then
\begin{equation}
\label{last3.5}
\lambda_{n, {\mathcal{N}}}[T_{\epsilon }(\Omega_1)]\le \lambda_{n, {\mathcal{N}}}[ \Omega_3](1+d_nA\epsilon ),
\end{equation}
where $d_n$ is defined by (\ref{patchesconst}). By inequality (\ref{last3.5}) and Lemma~\ref{unifbound}, it  follows
that there exist $c_n, \epsilon_n> 0$  such that
\begin{equation}
\label{last4}
\lambda_{n, {\mathcal{N}}}[T_{\epsilon }(\Omega_1)]\le \lambda_{n, {\mathcal{N}}}[ \Omega_3]+c_n\epsilon
\end{equation}
if $0<\epsilon <\epsilon_n$. On the other hand, by Lemma \ref{unifbound},  Corollary~\ref{diffeocorol}, and by inequalities  (\ref{buda1b}), (\ref{buda1,5b})
it follows that there exist $c_n$, $\epsilon_n>0$ such that
\begin{equation}
\label{last5}
| \lambda_{n, {\mathcal{N}}}[T_{\epsilon }[\Omega_1]]-\lambda_{n, {\mathcal{N}}}[\Omega_1]  |\le c_n \epsilon
\end{equation}
if $0<\epsilon <\epsilon_n$.
Thus by (\ref{last4}), (\ref{last5}) it  follows
that there exist $c_n, \epsilon_n> 0$  such that
\begin{equation}
\label{last6}
\lambda_{n, {\mathcal{N}}}[\Omega_1]\le \lambda_n[\Omega_3] +c_n\epsilon
\end{equation}
if $0<\epsilon <\epsilon_n $.
By Lemma \ref{semi} applied to the couple of open sets $\Omega_1,\Omega_3 $ it follows that there exist $c_n,\epsilon_n>0$ such that
\begin{equation}
\label{last6.5}
\lambda_{n, {\mathcal{N}}}[\Omega_3]\le \lambda_{n, {\mathcal{N}}}[\Omega_1] +c_n\epsilon
\end{equation}
if $0<\epsilon <\epsilon_n$. Thus, by (\ref{last6}), (\ref{last6.5}) it follows that
\begin{equation}
\label{last7}
|\lambda _{n, {\mathcal{N}}}[\Omega_1]-\lambda_{n, {\mathcal{N}}}[\Omega_3]  |\le c_n\epsilon
\end{equation}
if $0<\epsilon <\epsilon_n$. Clearly inequality (\ref{last7}) holds also with $\Omega_2$ replacing $\Omega_1$: it is simply enough to interchange the role of $\Omega_1$ and $\Omega_2$ from the beginning this proof. Thus
\begin{equation}
\label{last8}
|\lambda _{n, {\mathcal{N}}}[\Omega_2]-\lambda_{n, {\mathcal{N}}}[\Omega_3]  |\le c_n\epsilon
\end{equation}
if $0<\epsilon < \epsilon_n$.  By (\ref{last7}), (\ref{last8})
we finally deduce  that  for each $n\in {\mathbb{N}}$ there exist $c_n, \epsilon_n >0$  such that
\begin{equation}
\label{neuthm1bis}
| \lambda_{n, {\mathcal{N}}}[\Omega_1]-\lambda_{n  , {\mathcal{N}}}[\Omega_2] |\le c_n \epsilon ,
\end{equation}
for all $0<\epsilon <\epsilon_n$ and for all $\Omega_1, \Omega_2\in C({\mathcal{A}})$ satisfying $d_{\mathcal{A}}(\Omega_1,\Omega_2) <\epsilon $. Finally, by arguing as in the last lines of the proof of Theorem \ref{dirthm} we deduce the validity of (\ref{neuthm1}).
\hfill $\Box$\\

As for Dirichlet boundary conditions we have a version of Theorem \ref{neuthm} in terms of $\epsilon $-neighborhoods with respect to the atlas distance.

\begin{corol}
\label{lastcorANEU}
Let ${\mathcal{A}}$ be an atlas in ${\mathbb{R}}^N$. Let $m\in {\mathbb{N}}$,  $L, \theta >0$ and, for all $\alpha ,\beta \in {\mathbb{N}}_0^N$, with $|\alpha |=|\beta |=m$ let $A_{\alpha \beta }\in C^{0,1}(\cup_{j=1}^sV_j)$  satisfy $A_{\alpha \beta }=A_{ \beta \alpha }$, $\| A_{\alpha\beta} \|_{ C^{0,1}(\cup_{j=1}^sV_j)}\le L $ and condition (\ref{elp}).

 Then for each $n\in {\mathbb{N}}$ there exist $c_n, \epsilon_n >0$ depending only on
$n, N,   {\mathcal{A}}, m, L, \theta  $   such that
\begin{equation}
\label{maincordavstyAneu}
| \lambda_{n, {\mathcal{N}}}[\Omega_1]-\lambda_{n, {\mathcal{N}} }[\Omega_2] |\le c_n  \epsilon ,
\end{equation}
for all $0<\epsilon <\epsilon_n$ and
 for all $\Omega_1, \Omega_2\in C({\mathcal{A}})$ satisfying (\ref{maininclA}) or (\ref{maininclbisA}).
\end{corol}

{\bf Proof. } Inequality (\ref{maincordavstyAneu}) follows by inequality (\ref{neuthm1}) and  inequality (\ref{emibisA}). \hfill $\Box$\\

\section{Estimates via the lower Hausdorff-Pompeiu deviation}
\label{section7}

If $C\subset {\mathbb{R}}^N$ and $x\in {\mathbb{R}}^N$ we denote by $d(x,C)$ the euclidean distance of $x$ to $C$.

\begin{defn}\label{pompeiu}
Let $A,B \subset {\mathbb{R}}^N$. We define the
 lower Hausdorff-Pompeiu deviation of $A$ from $B$ by
\begin{equation}
\label{lowpo}
{\mathit{d}}_{{\mathcal{H}}{\mathcal{P}}} (A,B)=\min \biggl\{ \sup_{x\in A}d(x, B),\, \sup_{x\in B}d(x, A)\biggr\}.
\end{equation}
\end{defn}

If the minimum in (\ref{lowpo}) is replaced by the maximum, then the right-hand side becomes the usual Hausdorff-Pompeiu distance
${\mathit{d}}^{{\mathcal{H}}{\mathcal{P}}} (A,B)$
of $A$ and $B$. Note that in contrast to the Hausdorff-Pompeiu distance ${\mathit{d}}^{{\mathcal{H}}{\mathcal{P}}} $ which satisfies the triangle inequality and defines a distance on the family of closed sets, the lower Hausdorff-Pompeiu deviation  ${\mathit{d}}_{{\mathcal{H}}{\mathcal{P}}} $ is not a distance or a quasi-distance. Indeed, it suffices to notice that ${\mathit{d}}_{{\mathcal{H}}{\mathcal{P}}} (A,B)=0$ if and only if $A \subset \bar B$ or $B\subset \bar A$: thus if
$A \not\subset \bar B$ and $B \not\subset \bar A$ then ${\mathit{d}}_{{\mathcal{H}}{\mathcal{P}}} (A,B)>0 $ but
${\mathit{d}}_{{\mathcal{H}}{\mathcal{P}}} (A,A\cup B)+{\mathit{d}}_{{\mathcal{H}}{\mathcal{P}}} (A\cup B , A)=0$.

 In this section we aim at proving an estimate for  the variation of the eigenvalues  via  the lower Hausdorff-Pompeiu deviation
  of the boundaries of the open sets.

We now introduce a class of open sets for which  we can estimate the atlas distance $d_{{\mathcal{A}}}$  via the lower Hausdorff-Pompeiu deviation of the
boundaries.

\begin{defn}\label{omegafam}Let ${\mathcal{A}}$ be an atlas in ${\mathbb{R}}^N$.
Let $\omega :[0,\infty [\to [0,\infty [ $ be a continuous non-decreasing function
such that $\omega (0)=0$ and, for some $k>0$, $\omega (t)\geq kt$ for all $0\le t\le 1$.

Let $M>0$.
We denote by  $C_M^{\omega (\cdot )}( {\mathcal{A}}   )$ the family of all open sets  $\Omega$ in ${\mathbb{R}}^N$ belonging to
$C( {\mathcal{A}}   )$ and such that all the functions $g_j$ in Definition \ref{class} $(iii)$ satisfy the condition
\begin{equation}
\label{omegafamcond}
|g_j(\bar x) -g_j(\bar y)|\le M \omega (| \bar x -\bar y| ),
\end{equation}
for all $\bar x, \bar y \in {\overline{W}}_j$.

We also say that an open set is of class $C^{\omega (\cdot )}$ if there exists an atlas ${\mathcal{A}}$ and $M>0$ such that $\Omega \in C^{\omega (\cdot )}_M({\mathcal{A}} )$.

\end{defn}

\begin{lem}
\label{omega} Let  $\omega :[0,\infty [\to [0,\infty [$
be a continuous non-decreasing function
such that $\omega (0)=0$ and, for some $k>0$, $\omega (t)\geq kt$ for all $t\geq 0$.

Let $W$ be an open set in $\mathbb{R}^{N-1}$.
Let $M>0$ and
$g$ be a
function of ${\overline{W}}$ to $\mathbb{R}$ such
that
$$
|g(\bar x)-g(\bar y)|\le M\omega (|\bar x -\bar y|),
$$
for all $\bar x ,\, \bar y\in \overline{W}$. Then
\begin{equation}
\label{omega2}
 |g(\bar x)-x_N|\le (M+k^{-1})
 \omega \left( d\left((\bar x,x_N),{\rm
 Graph(g)}\right)\right),
\end{equation}
for all $\bar x\in\overline{ W}$ and $x_N\in\mathbb{R}$.

\end{lem}

{\bf Proof.} For all $\bar x, \bar y\in W$
\begin{eqnarray}\lefteqn{
|g(\bar x)-x_N|\le |g(\bar x)-g(\bar y)|+|g(\bar y)-x_N|}\nonumber \\  & & \qquad\qquad\qquad
\le M\omega (|\bar x-\bar y|)+k^{-1}\omega (|g(\bar y)-x_N|)\nonumber \\ & & \qquad\qquad\qquad
\le (m+k^{-1})\omega ( | (\bar x, x_N)-(\bar y, g(\bar y))  |    ),
\end{eqnarray}
hence by the continuity of $\omega $
\begin{eqnarray}
\lefteqn{
|g(\bar x)-x_N|\le  (m+k^{-1})\inf_{\bar y\in \overline{W}}\omega ( | (\bar x, x_N)-(\bar y, g(\bar y))  |              }\nonumber \\  & & \qquad\qquad\qquad
\le      (m+k^{-1})\omega ( \inf_{\bar y\in \overline{W}}| (\bar x, x_N)-(\bar y, g(\bar y))  |                          \nonumber \\ & & \qquad\qquad\qquad
\le (M+k^{-1})
 \omega \left( d\left((\bar x,x_N),{\rm
 Graph(g)}\right)\right).
\end{eqnarray}
\hfill $\Box$\\
\vspace{12pt}

\begin{lem}
\label{weighted-1}

Let ${\mathcal{A}}$  be an atlas in ${\mathbb{R}}^N$.   Let $\omega :[0,\infty [\to [0,\infty [$
be a continuous increasing function
satisfying $\omega (0)=0$ and, for some $k>0$, $\omega (t)\geq kt$ for all $0\le t\le 1$. Let $M>0$.

Then there exists $c >0$ depending only on $N, {\mathcal{A}},
 \omega , M$  such that
\begin{equation}
\label{weighted1} d_j(x, \partial\Omega)\le c\, \omega
(d (x, \partial\Omega ) ),
\end{equation}
for all open sets $\Omega \in C^{\omega (\cdot )}_M({\mathcal{A}} )$,   for all $j= 1,
\dots, s$ and for all $x\in (V_j)_{\frac{\rho}{2}}$.
\end{lem}

{\bf Proof.}
Let $\tilde \omega$ be the function of $[0, \infty [$ to itself defined by
$\tilde \omega (t) =\omega (t) $ for all $0\le t \le 1$ and $\tilde \omega (t)= t+\omega (1) -1$ for all
$t>1$. Clearly $\tilde \omega $ is continuous and non-decreasing,   and $\tilde \omega (t)\geq \tilde k t$ for all $t\geq 0$ where $\tilde k= \min \{k,1, \omega (1)\}$; moreover
\begin{equation}
\label{tildeomega}  \min \left\{ 1, \frac{\omega (1)}{\omega (A)}     \right\}\omega (a)\le
\tilde \omega (a) \le \max\left\{ 1, \frac{\tilde \omega (A)}{\tilde \omega (1)}  \right\}\omega (a),
\end{equation}
for all $A>0$ and for all $0\le a\le A$.
By the first inequality in (\ref{tildeomega}) it follows that  $\Omega\in {\rm C}^{\tilde \omega (\cdot ) }_{\tilde M} ( {\mathcal{A}} )$ where $\tilde M =  \max\left\{ 1, \frac{ \omega (D)}{\omega (1)}  \right\} M$ and $D$ is the diameter of $\cup_{j=1}^{s} V_j$.
Then by  Lemma \ref{omega} applied to each function $g_j$ describing the boundary of $\Omega $ as in Definition \ref{class} $(iii)$ and by the second inequality in (\ref{tildeomega}) it follows that for each $j=1, \dots s$ and for all $(\bar y, y_N )\in r_j(V_j)$
\begin{eqnarray}
\label{tildeomega0}
\lefteqn{
|g_j(\bar y)-y_N|\le  (\tilde M+\tilde k^{-1})
\, \tilde \omega \left( d\left((\bar y,y_N),{\rm
 Graph}(g_j)\right)\right) }\nonumber \\ & & \, \le  \max\left\{ 1, \frac{\tilde \omega ( D)}{\tilde \omega (1)}  \right\}
 (\tilde M+\tilde k^{-1})\,
  \omega \left( d\left((\bar y,y_N),{\rm
 Graph}(g_j)\right)\right) .
\end{eqnarray}
Observe that if $y\in r_j((V_j)_{\frac{\rho}{2}})$ and  $ d\left((\bar y,y_N),{\rm
 Graph}(g_j)\right)< \frac{\rho }{2}$ then  $d(r_j^{(-1)}(y),\partial \Omega ) $ equals $d\left((\bar y,y_N),{\rm
 Graph}(g_j)\right)$;  if $y\in r_j((V_j)_{\frac{\rho}{2}})$ and  $ d\left((\bar y,y_N),{\rm
 Graph}(g_j)\right)\geq \frac{\rho }{2}$ then  $ d\left(r_j^{(-1)}((\bar y,y_N)), \partial \Omega \right)\geq \frac{\rho }{2}$. Hence
\begin{equation}
\label{tildeomega1}
 \omega \left( d\left((\bar y,y_N),{\rm
 Graph}(g_j)\right)\right)\le \frac{ \omega \left( D\right)}{\omega ( \frac{\rho}{2}) }     \omega \left( d\left(r_j^{(-1)}((\bar y,y_N)),\partial \Omega \right)\right).
\end{equation}
So by (\ref{tildeomega0}) and (\ref{tildeomega1})
it follows that if $y\in r_j((V_j)_{\frac{\rho}{2}})$ then
$$
|g_j(\bar y) -y_N|\le
c
 \omega \left( d(r_j^{(-1)}(y),\partial \Omega )\right),
$$
where $c=\max\left\{ 1, \tilde \omega ( D)/\tilde \omega (1)  \right\}
 (\tilde M+\tilde k^{-1})   \omega \left( D\right)/\omega ( \rho /2) $. Hence, for $x\in (V_j)_{\frac{\rho}{2}}$, by (\ref{dj})
 $$
 d_j(x, \partial\Omega )=|g_j(\overline{r_j(x)})-(r_j(x))_N|\le c\omega (d(x, \partial\Omega )).
 $$
\hfill $\Box$\\

\begin{lem}
\label{hausdorff00}
Let ${\mathcal{A}} =(    \rho  , s,s', \{ (V_j)_{j=1}^s, \{r_j\}_{j=1}^{s} )   $  be an atlas in ${\mathbb{R}}^N$. Let $\tilde {\mathcal{A}}=(    \rho /2 , s,s', \{ (V_j)_{\rho /2}\}_{j=1}^s, \{r_j\}_{j=1}^{s} )   $.  Let $\omega :[0,\infty [\to [0,\infty [$
be a continuous non-decreasing  function
satisfying $\omega (0)=0$ and, for some $k>0$, $\omega (t)\geq kt$ for all $0\le t\le 1$. Let $M>0$.

Then there exists $c >0$ depending only on $N, {\mathcal{A}} ,
 \omega ,  M$  such that
\begin{equation}
\label{hausdorff}
d^{{\mathcal{H}}   {\mathcal{P}}}(\partial \Omega_1,\partial \Omega_2)\le
d_{\tilde {\mathcal{A}}}(\Omega_1, \Omega_2) \le c\, \omega (d_{{\mathcal{H}}   {\mathcal{P}}}(\partial \Omega_1,\partial \Omega_2)),
\end{equation}
for all opens sets $\Omega_1, \Omega_2 \in C^{\omega (\cdot ) }_M({\mathcal{A}})$.
\end{lem}

{\bf Proof. }  For each $x\in \partial \Omega_1$ there exists $y\in \partial \Omega_2$ such that $|x-y|\le d_{\tilde {\mathcal{A}}}(\Omega_1, \Omega_2)$: indeed, if $r_j(x)=(\bar x , x_N)$ for some $j=1, \dots , s'$ it is sufficient to consider  $y\in \partial \Omega_2$ such that
$\overline {r_j(y)}=\bar x$.  It follows that $d(x, \partial \Omega_2)\le d_{\tilde {\mathcal{A}}}(\Omega_1, \Omega_2)$ for all $x\in \partial \Omega_1$. In the same way, $d(x, \partial \Omega_1)\le d_{\tilde {\mathcal{A}}}(\Omega_1, \Omega_2)$ for all $x\in \partial \Omega_2$. Thus,  the first inequality in (\ref{hausdorff}) follows.
The second inequality in (\ref{hausdorff}) immediately follows by (\ref{weighted1}), by the continuity of $\omega $, by property (\ref{atl}) and by Definition~\ref{pompeiu}.\hfill $\Box$\\

\begin{thm}\label{maincor}
Let ${\mathcal{A}}$ be an atlas in ${\mathbb{R}}^N$. Let $m\in {\mathbb{N}}$,  $L, M,\theta >0$ and, for all $\alpha ,\beta \in {\mathbb{N}}_0^N$ with $|\alpha |=|\beta |=m$, let $A_{\alpha \beta }\in C^{0,1}(\cup_{j=1}^sV_j)$  satisfy $A_{\alpha \beta }=A_{ \beta \alpha }$, $\| A_{\alpha\beta} \|_{ C^{0,1}(\cup_{j=1}^sV_j)}\le L $ and condition (\ref{elp}). Let $\omega :[0,\infty [\to [0,\infty [$
be a continuous non-decreasing function
satisfying $\omega (0)=0$ and, for some $k>0$, $\omega (t)\geq kt$ for all $0\le t\le 1$.

 Then for each $n\in {\mathbb{N}}$ there exist $c_n, \epsilon_n >0$ depending only on
$n, N,   {\mathcal{A}}, m,L, M,$ $ \theta , \omega $   such that for both  Dirichlet and  Neumann boundary conditions
\begin{equation}
\label{neuthm1ha}
| \lambda_{n}[\Omega_1]-\lambda_{n }[\Omega_2] |\le c_n \omega ( d_{{\mathcal{H}}   {\mathcal{P}}}(\partial \Omega_1,\partial \Omega_2)) ,
\end{equation}
 for all $\Omega_1, \Omega_2\in C_M^{\omega (\cdot)}({\mathcal{A}})$ satisfying $d_{{\mathcal{H}} {\mathcal{P}}}(\partial\Omega_1,\partial \Omega_2) <\epsilon_n $.
\end{thm}

{\bf Proof.} Observe  that if $\Omega_1,\Omega_2\in C({\mathcal{A}})$ then also
$\Omega_1,\Omega_2\in C(\tilde {\mathcal{A}})$ where $\tilde {\mathcal{A}}=(    \rho /2 , s,s', \{ (V_j)_{\rho /2}\}_{j=1}^s, \{r_j\}_{j=1}^{s} )   $. Thus by inequalities (\ref{dirthm1}),
(\ref{neuthm1}) applied to  $\Omega_1, \Omega_2$ as open sets in  $C(\tilde {\mathcal{A}})$ and by inequality (\ref{hausdorff})  we deduce the validity of (\ref{neuthm1ha}). \hfill $\Box$\\

Recall that for any $\Omega$ we set
$$
\Omega ^{\epsilon }=\{x\in {\mathbb{R}}^N:\ d(x, \Omega )<\epsilon   \},
$$
and
$$
\Omega _{\epsilon }=\{x\in \Omega :\ d(x, \partial \Omega )>\epsilon   \}.
$$

\begin{lem} If $\Omega_1$ and $\Omega_2$ are
two open sets satisfying the inclusions
\begin{equation} \label{mainincl}
(\Omega_1)_{\epsilon }\subset \Omega_2 \subset (\Omega_1)^{\epsilon }
\end{equation}
or
\begin{equation} \label{maininclbis}
(\Omega_2)_{\epsilon }\subset \Omega_1 \subset (\Omega_2)^{\epsilon },
\end{equation}
then
\begin{equation}
\label{emibis}
d_{{\mathcal{H}}{\mathcal{P}}}(\partial \Omega_2,\partial \Omega_1)\le \epsilon .
\end{equation}
\end{lem}

{\bf Proof.}  Similarly to the proof of Lemma~\ref{cinquesedici} inclusions (\ref{mainincl}) and (\ref{maininclbis}) imply that
$\sup_{x\in \partial \Omega_2}d(x,\partial\Omega_1)\le \epsilon $, $\sup_{x\in \partial \Omega_1}d(x,\partial\Omega_2)\le \epsilon $
respectively. Hence if (\ref{mainincl}) or (\ref{maininclbis}) is satisfied then at least one of these inequalities is satisfied
which implies (\ref{emibis}).\hfill $\Box$\\

Observe that if $\Omega_1$ and $\Omega_2$ are
two open sets satisfying inclusion (\ref{mainincl})  then it may happen that they do not satisfy inclusion (\ref{maininclbis}), and
\begin{equation}
\label{emibisbis}
\sup_{x\in \partial \Omega_1}d(x, \partial \Omega_2 )> \epsilon ,
\end{equation}
see Examples \ref{ex1}, \ref{ex2},  in Appendix.

\begin{corol}
\label{lastcor} Under the same assumptions of Theorem~\ref{maincor},
for each $n\in {\mathbb{N}}$ there exist $c_n, \epsilon_n >0$ depending only on
$n, N,   {\mathcal{A}}, m, L, M,$ $ \theta , \omega $   such that for both  Dirichlet and  Neumann boundary conditions
\begin{equation}
\label{maincordavsty}
| \lambda_{n}[\Omega_1]-\lambda_{n }[\Omega_2] |\le c_n \omega (  \epsilon ),
\end{equation}
for all $0<\epsilon <\epsilon_n$ and
 for all $\Omega_1, \Omega_2\in C_M^{\omega (\cdot )}({\mathcal{A}})$ satisfying (\ref{mainincl}) or (\ref{maininclbis}).
\end{corol}

{\bf Proof. } Inequality (\ref{maincordavsty}) follows by inequalities (\ref{neuthm1ha}) and (\ref{emibis}).\hfill $\Box$\\

\section{Appendix }

\subsection{On the atlas distance}

Given an atlas  ${\mathcal{A}}$  in ${\mathbb{R}}^N$, it is immediate to prove that the function $d_{{\mathcal{A}}}$ of
$ C({\mathcal{A}})\times C({\mathcal{A}})$ to ${\mathbb{R}}$ which takes $(\Omega_1,\Omega_2)$ to $d_{{\mathcal{A}}}(\Omega_1,\Omega_2)$
for all  $(\Omega_1,\Omega_2)\in C({\mathcal{A}})\times C({\mathcal{A}})$, is a metric on the set $C({\mathcal{A}})$.

\begin{lem}
\label{convergenza}
Let ${\mathcal{A}}$ be an atlas in ${\mathbb{R}}^N$. Let  $\Omega _n, \ n\in {\mathbb{N}}$,  be a sequence in
$C({\mathcal{A}})$. For each $n\in {\mathbb{N}}$ let $g_{jn}$, $j=1, \dots ,s$, be the functions describing the boundary of $\Omega_n$
as in Definition \ref{class} (iii). Then the sequence $\Omega _n, \ n\in {\mathbb{N}}$, is convergent in $(C({\mathcal{A}}), d_{{\mathcal{A}}})$ if and only if for all $j=1, \dots , s$ the sequences $g_{jn},\ n\in {\mathbb{N}}$, are uniformly convergent on $\overline{ W}_j$. Moreover,
if $g_{jn} $ converge uniformly to $g_j$ on $\overline {W}_j$ for all $j=1, \dots , s$ then $\Omega_n $ converges in $(C({\mathcal{A}}), d_{{\mathcal{A}}}  )$ to the open set $\Omega \in C({\mathcal{A}})$ whose boundary is described
by the functions $g_j$ as in Definition \ref{class} (iii).
\end{lem}

{\bf Proof. } It is enough to prove that  if the sequences $g_{jn},\ n\in {\mathbb{N}} $, converge  to $g_j$ uniformly on $\overline {W}_j$ for all $j=1, \dots , s$ then the sequence $\Omega_n $, $n\in {\mathbb{N}}$, converges in $(C({\mathcal{A}}), d_{{\mathcal{A}}})$ to the open set $\Omega\in C({\mathcal{A}}) $ whose boundary is described
by the functions $g_j$ as in Definition \ref{class} $(iii)$ (the rest is obvious). We divide the proof into two steps.

{\it Step 1.}
 We  prove that if  $x\in V_h\cap V_k$ for $h\ne k$ and  $r_h(x)=  (\bar x , g_h(\bar x)) $ for some $\bar x\in W_h$ then there exists
 $\bar y\in W_k$ such that  $r_k(x)=(\bar y, g_k(\bar y) )$.
 Observe that $x=\lim_{n\to \infty }r_h^{(-1)}(\bar x,  g_{hn}(\bar x))$ and there exists $\tilde n\in {\mathbb{N}}$ such that  $r_h^{(-1)}(\bar x,  g_{hn}(\bar x))\in V_h\cap V_k$ for all $n\geq \tilde n$. For each $n\geq \tilde n$ there exists $\bar y_n\in W_k$ such
that $r_k (r_h^{(-1)}(\bar x,  g_{hn}(\bar x)))=(\bar y_n, g_{kn}( \bar y_n ))$. Clearly $\lim_{n\to \infty} r_k (r_h^{(-1)}(\bar x,  g_{hn}(\bar x)))=r_k(x) $ hence  $\lim_{n\to\infty }(\bar y_n, g_{kn}( \bar y_n ))=r_k(x)$. By the uniform convergence of $g_{kn}$ to
$g_k $ on $W_k$ it follows that there exists $\bar y\in W_k$ such that $\lim_{n\to \infty }(\bar y_n, g_{kn}(\bar y ))=(\bar y, g_k(\bar y)) $. Thus
$\lim_{n\to\infty } r_k( r_h^{(-1)}(\bar x,  g_{hn}(\bar x)))=(\bar y, g_k(\bar y)) $ and
$x=r_k^{(-1)}(\bar y,  g_{k}(\bar y))$ as required.

{\it Step 2.}
We  prove  that if $x\in V_h\cap V_k$ for $h\ne k$,  $r_h(x)=  (\bar x , x_N )$ for some $\bar x\in W_h$ and $x_N<g_h(\bar x)$ then there exists
 $\bar y\in W_k$ such that
$r_k(x)=(\bar y, y_N)$ and $ y_N <g_k(\bar y)$. Indeed there exists $\hat n\in {\mathbb {N}}$ such that
$ x_ N <g_{hn}(\bar x ) $  for all $n\geq \hat n$. Thus $x\in V_h\cap V_k\cap  \Omega_n$, hence
$$
  (r_k(x) )_N<g_{kn }(\overline{ r_k(x)}),
$$
and by passing to the limit it follows that
$$
   (r_k(x) )_N\le g_{k }(\overline{ r_k(x)}).
$$
If  $ (r_k(x) )_N  =g_{k }(\overline{ r_k(x)}) $, then by {\it Step 1} there exists $\bar z\in W_h$ such that $r_h(x)=(\bar z , g_h(\bar z )  ) $ which implies $\bar z=\bar x$ and $g_h(\bar z)=x_N$
which contradicts the assumption that $x_N<g_h(\bar x)$.
Thus we have proved that
$$
 (r_k(x) )_N<g_{k }(\overline{ r_k(x)}).
$$
In other words, $r_k(x)= (\bar y ,y_N )$ where $\bar y= \overline{ r_k(x)}$, $y_N=(r_k(x) )_N$, and $y_N<g_k(\bar y)$ as required.

By {\it Steps 1,2}\quad it follows that the set
$$
\Omega =\bigcup_{j=1}^s r_j^{(-1)}\left( \left\{(\bar x, x_N  ):\ \bar x\in W_j,\ a_{Nj}< x_N< g_j(\bar x)\right\} \right)
$$
is such that
$$
 r_j(   \Omega \cap V_j) =  \left\{(\bar x, x_N):\ \bar x\in W_j,\ a_{Nj} < x_N< g_j(\bar x) \right\}.
$$
Thus $\Omega \in C({\mathcal{A}})$. Obviously $\lim_{n\to \infty }d_{{\mathcal{A}}}(\Omega_n, \Omega )=0$. \hfill $\Box$

\begin{thm}Let ${\mathcal{A}}$ be an atlas in ${\mathbb{R}}^N$. Then $(C({\mathcal{A}}), d_{{\mathcal{A}}})$ is a complete
metric space. Moreover, for each function $\omega$ satisfying the assumptions of Definition~\ref{omegafam} and for each $M>0$,  $C_M^{\omega (\cdot )} ({\mathcal{A}})$ is a compact set in $(C({\mathcal{A}}), d_{{\mathcal{A}}})$.
\end{thm}

{\bf Proof.}  The completeness of  $(C({\mathcal{A}}), d_{{\mathcal{A}}})$ and the closedness of the set $C_M^{\omega (\cdot )} ({\mathcal{A}})$
follows directly  by Lemma~\ref{convergenza} (in the second case one should take into account that condition  (\ref{omegafamcond}) with fixed
$\omega$ and $M$ allows passing to the limit).

By the definition of  $C_M^{\omega (\cdot )} ({\mathcal{A}})$ it follows that for all $j=1,\dots , s$
the sets $G_j=\{g_j[\Omega ]\}_{\Omega \in C_M^{\omega (\cdot )} ({\mathcal{A}})}$ of functions $g_j[\Omega ]$ entering Definition~\ref{class}, which are defined on the bounded cuboids $\overline{W}_j$, are bounded in the sup-norm and equicontinuous due to condition (\ref{omegafamcond}) where $\omega $ and $M$ are the same for all $\Omega \in C_M^{\omega (\cdot )} ({\mathcal{A}}) $. By the Ascoli-Arzel\`{a} Theorem the sets $G_j$, being closed,
are compact with respect to the sup-norm.

Let $\{\Omega_n\}_{n\in {\mathbb{N}}}$ be a sequence in $C_M^{\omega (\cdot )} ({\mathcal{A}})$. Since the sets $G_j$ are compact it follows that,
possibly considering a subsequence, $\{ g_{jn}\}_{n\in {\mathbb{N}}}$, where $ g_{jn}=g_{j}[\Omega_n ]$, converges uniformly on $\overline {W}_j$ to
some continuous functions $g_j$, $j=1, \dots , s$. By Lemma~\ref{convergenza} the sequence $\{\Omega_n\}_{n\in {\mathbb{N}}}$ converges
in  $(C({\mathcal{A}}), d_{{\mathcal{A}}})$ to the open set $\Omega $ defined by the functions $g_j$, $j=1, \dots , s$. Therefore the set
$C_M^{\omega (\cdot )} ({\mathcal{A}})$ is relatively  compact in $(C({\mathcal{A}}), d_{{\mathcal{A}}})$ and, being closed, is compact.\hfill $\Box$\\

\subsection{Comparison of atlas distance, Hausdorff-Pompeiu distance, and lower Hausdorff-Pompeiu deviation}

Observe that

\begin{equation}
\label{dist}
d_{{\mathcal{H}}{\mathcal{P}}} (\partial \Omega_1,\partial \Omega_2) \le
d^{{\mathcal{H}}{\mathcal{P}}} (\partial \Omega_1, \partial \Omega_2)\le   d_{\mathcal{A}}( \Omega_1, \Omega_2 )
,
\end{equation}
for all $\Omega_1,\Omega_2\in C({\mathcal{A}})$. (The first inequality is trivial, the second can be proved by using the same argument used in the proof of Lemma \ref{hausdorff00}.)

The following examples show that
$d_{{\mathcal{H}}{\mathcal{P}}} (\partial \Omega_1,\partial \Omega_2)$ can be much smaller than  $d^{{\mathcal{H}}{\mathcal{P}}} (\partial \Omega_1, \partial \Omega_2)$, and $d^{{\mathcal{H}}{\mathcal{P}}} (\partial \Omega_1, \partial \Omega_2)$ can be much smaller than
$d_{\mathcal{A}}( \Omega_1, \Omega_2 )$.

\begin{example}\label{ex1} Let $N=2$, $c>\sqrt{3}$ and $0<\epsilon < 1/\sqrt{3}$. Let
$$
\Omega_1 =\left\{(x_1, x_2)\in {\mathbb{R}}^2:\ c|x_2| < x_1 <c  \right\}
$$
and $\Omega_2=(\Omega_1)_{\epsilon }$. Then $\Omega_1, \Omega _2$ satisfy inclusion (\ref{mainincl}) but not (\ref{maininclbis}).

Moreover, the lower Hausdorff-Pompeiu deviation of the boundaries can be much smaller than their usual Hausdorff-Pompeiu distance, because
$$
d_{{\mathcal{H}}{\mathcal{P}}} (\partial \Omega_1, \partial\Omega_2)=\epsilon\ \ {\rm and }\ \ d^{{\mathcal{H}}{\mathcal{P}}} (\partial \Omega_1, \partial \Omega_2) = \epsilon \sqrt{c^2+1}.
$$
\end{example}

\begin{example}\label{ex2}
Let a function $\psi : [0, \infty [\to [0, \infty [$ be such that $\psi (0)=\psi'(0)=0$ and $\psi '(t)>0$ for all $t>0$. Let $\omega : [0, \infty [\to [0, \infty [$ be the inverse function of $\psi$. Let $N=2$ and
$$
\Omega_1 =\{ (x_1 ,x_2)\in {\mathbb{R}}^2:\ \omega (|x_2|)<x_1 <\omega (1)  \}.
$$
Let $P=(x_1 ,0)$ be a point with   $x_1>0$ sufficiently small so that $d(P, \partial \Omega_1 )= d(P, \{(t, \psi (t)):\ 0\le t\le \omega (1)\})=|P-Q|$ for some point $Q=(\xi , \psi (\xi ))$ with $0<\xi < \omega (1)$.
We set $d(P, \partial \Omega_1 )=\epsilon $. By elementary considerations it follows that
$$
x_1= \xi + \psi (\xi )\psi '(\xi )\ \ {\rm and }\ \ \epsilon =\psi (\xi )\sqrt {1+\psi '(\xi )^2}.
$$
This implies that $x_1\sim \xi$ and $\epsilon \sim \psi (\xi ) $ as $\xi \to 0^{+}$, hence $x_1\sim \omega (\epsilon ) $ as $\epsilon \to 0^{+}$.

Now let $\Omega_2= (\Omega _1)_{\epsilon }$.
Clearly $\Omega_1, \Omega _2$ satisfy inclusion (\ref{mainincl})
and
$$
  \sup_{x\in \partial \Omega_2}d(x, \partial \Omega_1)=\epsilon .
$$
However, since   $P\in \partial \Omega_2$, we have
$$
 \sup_{x\in \partial \Omega_1}d(x, \partial \Omega_2)\geq x_1 \sim \omega (\epsilon )
$$
as $\epsilon \to 0^{+}$,
hence $\Omega_1, \Omega_2$ cannot satisfy  (\ref{maininclbis}) for small values of $\epsilon$ because $\lim_{\epsilon \to 0^{+}} \omega (\epsilon )/\epsilon =\infty $. Moreover, there exist $c_1, c_2>0$ such that for all sufficiently small $\epsilon >0$
$$c_1\omega( d_{{\mathcal{H}}{\mathcal{P}}} (\partial \Omega_1,\partial \Omega_2) )\le
d^{{\mathcal{H}}{\mathcal{P}}} (\partial \Omega_1, \partial \Omega_2)\le c_2 \omega( d_{{\mathcal{H}}{\mathcal{P}}} (\partial \Omega_1,\partial \Omega_2) ).
$$
(The second inequality follows by (\ref{hausdorff}).)
 This means in particular that  the usual Hausdorff-Pompeiu distance  $ d^{{\mathcal{H}}{\mathcal{P}}} $  between the boundaries may  tend to zero  arbitrarily slower than  their lower Hausdorff-Pompeiu deviation $d_{{\mathcal{H}}{\mathcal{P}}} $.

\end{example}

\begin{example}
\label{atlasvshaus}
Let $N= 2$. Let ${\mathcal{A}}= (  \rho , s,s', \{V_j\}_{j=1}^s, \{r_j\}_{j=1}^{s} ) $  be an atlas in ${\mathbb{R}}^2$ with $V_1=]-2,2[\times ]-2,2[$. Let $\omega :[0,\infty [\to [0,\infty [ $ be a continuous increasing function
such that $\omega (0)=0$ and, for some $k>0$, $\omega (t)\geq kt$ for all $0\le t\le 1$. Assume also that $\omega (1)=1$ and
that, for some $M>0$, $|\omega (x)-\omega (y)|\le M\omega (|x-y|)$ for all $0\le x, y \le 1$. Let $0<\epsilon <1/2$.
Let $\Omega_1, \Omega_2\in C^{\omega (\cdot )}_M({\mathcal{A}})$ with
\begin{eqnarray}
r_1(\Omega_1\cap V_1)& =& \{( x_1, x_2):\ -2<x_1<2,\ -2<x_2<g_{11}( x_1) \}\nonumber \\
r_1(\Omega_2\cap V_1)& =& \{( x_1, x_2):\ -2<x_1<2,\ -2<x_2<g_{12}( x_1) \}\nonumber
\end{eqnarray}
where
\begin{eqnarray}
g_{11}(x_1)& =& \left\{\begin{array}{ll}  1-\omega (|x_1|)& {\rm if}\ |x_1|\le 1, \\  0 & {\rm if}\ 1< |x_1|< 2, \end{array}\right. \\
g_{12}(x_1)& =& \left\{\begin{array}{ll}g_{11}(x_1-\epsilon ) & {\rm if}\ -2+\epsilon \le x_1< 2, \\ 0& {\rm if}\ -2 <x_1 < -2+\epsilon \, , \end{array}\right.
\end{eqnarray}
and $\Omega_1\cap V_j=\Omega _2\cap V_j$ for all $2\le j \le s$. It is clear that
$ d^{{\mathcal{HP}}}(\partial \Omega_1,\partial \Omega_2)\le \epsilon $, and
$
d_{{\mathcal{A}}}(\Omega_1,\Omega_2)\geq g_{11}(0)-g_{12}(0)=\omega (\epsilon ).
$
Thus
\begin{equation}
\label{atlasvshaus1}
 \omega ( d^{{\mathcal{HP}}}(\partial \Omega_1,\partial \Omega_2) ) \le d_{{\mathcal{A}}}(\Omega_1,\Omega_2).
\end{equation}
\end{example}

{\bf Acknowledgements.}
This research was supported by the research project ``Problemi di stabilit\`{a} per operatori differenziali" of the University of Padova, Italy and
 partially supported by the grant of INTAS (project 05-1000008-8157). The first author was also supported by the grant of RFBR - Russian Foundation for Basic Research (project 05-01-01050).

{DIPARTIMENTO DI MATEMATICA PURA ED APPLICATA, UNIVERSIT\`{A} DEGLI STUDI DI PADOVA, VIA TRIESTE 63, 35121 PADOVA, ITALY}

\end{document}